\documentclass{article}

\usepackage[margin=1in]{geometry}
\usepackage{amsmath,amsthm,amssymb}
\usepackage{graphicx}
\usepackage{color}
\usepackage[dvipsnames]{xcolor}
\usepackage{subcaption}
\usepackage[export]{adjustbox}
\usepackage{url}

\graphicspath{ {./nlmc/} }

\begin{document}

\title{Nonlinear nonlocal multicontinua upscaling framework and its applications}
\author{
Wing T. Leung\thanks{ICES, University of Texas, Austin, TX, USA (\texttt{wleung@ices.utexas.edu})}
\and
 Eric T. Chung\thanks{Department of Mathematics, The Chinese University of Hong Kong, Shatin, New Territories, Hong Kong SAR, China (\texttt{tschung@math.cuhk.edu.hk}) }
\and
Yalchin Efendiev\thanks{Department of Mathematics \& Institute for Scientific Computation (ISC), Texas A\&M University,
College Station, Texas, USA (\texttt{efendiev@math.tamu.edu})}
\and
 Mary Wheeler\thanks{ICES, University of Texas, Austin, TX, USA (\texttt{mfw@ices.utexas.edu})}
}

\maketitle

\begin{abstract}

In this paper, we discuss multiscale methods for nonlinear problems. We use
recently developed multiscale concepts for linear problems and extend them
to nonlinear problems. The main idea of these approaches is to use
local constraints and solve problems in oversampled regions for
constructing macroscopic equations. These techniques are intended
for problems without scale separation and high contrast, which often
occur in applications. For linear
problems, the local solutions with constraints
 are used as basis functions. This technique is
called Constraint Energy Minimizing Generalized Multiscale Finite
Element Method (CEM-GMsFEM). GMsFEM identifies macroscopic
quantities based on rigorous analysis. In corresponding upscaling methods,
the multiscale basis functions are selected such that the degrees of
freedom have physical meanings, such as averages of the solution on
each continuum.

This paper extends the linear concepts to nonlinear problems,
where the local problems are nonlinear.
The main concept consists of: (1) identifying macroscopic quantities;
(2) constructing appropriate oversampled local  problems with
coarse-grid constraints; (3) formulating macroscopic equations.
We consider two types of approaches.
In the first approach, the solutions of local problems are used as basis
functions (in a linear fashion) to solve nonlinear problems.
This approach is simple to implement; however, it lacks the nonlinear
interpolation, which we present in our second approach. In this approach,
the local solutions are used as a nonlinear forward map from local
averages (constraints) of the solution in oversampling region.
This local fine-grid solution is further used to formulate the coarse-grid
problem. Both approaches are discussed on several examples and applied
to single-phase and two-phase flow problems, which are challenging
because of convection-dominated nature of the concentration equation.
The numerical results show that we can achieve  good accuracy using
our new concepts for these complex problems.


\end{abstract}

\section{Introduction}

As multiscale methods for linear equations are getting matured,
their extension to difficult
nonlinear problems remain challenging. Many nonlinear
problems have multiscale nature due to spatial and temporal scales.
For example, the dynamics of
multi-phase flow and transport in heterogeneous media
varies over multiple space and time scales.
In the past,
many well-known linear and
nonlinear upscaling tools have been developed
(e.g., \cite{ab05, egw10,  arbogast02, GMsFEM13, AdaptiveGMsFEM, brown2014multiscale, ElasticGMsFEM, ee03, abdul_yun, ohl12, fish2004space, fish2008mathematical, oz07, matache2002two, apwy07, henning2009heterogeneous, OnlineStokes, chung2017DGstokes,WaveGMsFEM, pwy02, Arbogast_PWY_07, MsDG, fish1997computational,oskay2007eigendeformation,yuan2009multiple}).
Single-phase upscaling methods include permeability upscaling
\cite{dur91, weh02, cdgw03,hn00}
and many multiscale techniques
\cite{ab05, GMsFEM13, AdaptiveGMsFEM, oz07, matache2002two, apwy07, henning2009heterogeneous}).
 Nonlinear
upscaling methods, e.g., known as pseudo-relative permeability approach
\cite{cdgw03,KB,BT},
computes nonlinear relative permeability functions based on single cell
two-phase flow computations.
It is known that these nonlinear approaches lack
robustness and they are processes dependent \cite{ed01,ed03}. To overcome
these difficulties, one needs a better understanding of nonlinear upscaling
methods, which is our goal.

Nonlinear upscaling methods can be traced back to nonlinear homogenization
\cite{pankov97, ep03d}.
The main idea of nonlinear homogenization is to formulate coarse-grid equations
based on nonlinear local problems formulated in each coarse block.
Some examples include p-Laplacian, pseudo-elliptic equations, and
parabolic equations \cite{ep03d}.
In these approaches, the local problems are
solved with periodic boundary conditions and some constraints on
averages of the solutions of gradients. In non-periodic cases, these
methods are extended by solving local problems in coarse block
subject to some boundary conditions. Because of
 nonlinearity, these local problems need to be solved for
all possible average values, which can make the computations expensive.
One can compute the upscaled fluxes on-the-fly using the values at previous
iteration or previous time. These approaches are limited to problems
without high contrast and scale separation. Our goal is to extend
these approaches to problems with high contrast and non-separable scales.
{\it The main novel components of our approach} is
introducing multiple macroscopic variables for each coarse-grid block,
formulating appropriate local constrained problems that determine
the downscaling map; formulating macroscopic equations.

In computational mechanics literature, there has been a great deal of research
dedicated to nonlinear upscaling methods, which include
generalized continuum theories (e.g., \cite{fafalis2012capability}),
computational continua
framework (e.g., \cite{fish2010computational}),
and other approaches.
Multiscale enrichment method based on the partition of unity (\cite{fish2005multiscale,fish2007multiscale}) combine the linear and
non-linear homogenization theory with the partition of unity method
to handle the problems with inseparable fine and coarse scales.
Computational continua (\cite{fish2010computational,fafalis2018computational}), which use nonlocal quadrature to couple the coarse scale system stated on a unions of some disjoint computational unit cells, are introduced for non-scale-separation heterogeneous media. In \cite{fish2012reduced,fish2015computational,fish2013practical}, the authors further enhance the method by combining the
computational continua with model reduction technique. Originally,
the computational continua approaches were
developed to overcome the limitations of higher-order
homogenization models and generalized continuum theories, namely,
the need for higher-order finite element continuity, additional degrees of
freedom, and nonclassical boundary conditions. Though there
are some similarities,
our proposed approaches differ from computational continua framework. Our
approach  explores the localization of high-contrast
(with respect to the coarse-mesh size)
features by introducing additional degrees of freedom (continua) and
using oversampled regions. In future, we plan to use some ingredients
of computational continua approaches in improving our approach and in making
it more applicable to computational materials.

To design the new upscaled model, we use the concept
of non-local multi-continuum (NLMC). The main idea of NLMC derives
from Constraint Energy Minimizing Generalized Multiscale Finite Element
Method (CEM-GMsFEM) \cite{chung2018constraint,chung2018fast,chung2018constraintmixed}.
CEM-GMsFEM identifies the degrees of freedom
which can not be localized and then construct multiscale basis functions
with support in the oversampled regions. However, the degrees
of freedom in CEM-GMsFEM, in general,
 do not have physical meanings since they are
coordinates in the multiscale space. In order to obtain an upscaled model,
we design multiscale basis functions such that the degrees of freedom
represent the average of the solutions. As a result, the coarse-grid model
is an equation for the solutions averaged over each continua. In this paper,
we extend this idea to nonlinear equations.

To extend the concept of NLMC to nonlinear equations, we
first identify macroscopic quantities for each coarse-grid block.
These variables are typically
 found via local spectral decomposition and
 represent the features that can not be localized
(similar to multicontina variables \cite{barenblatt1960basic,lee2001hierarchical, warren1963behavior, panasenko2018multicontinuum}).
Next, we consider
local problems formulated in the oversampled regions with constraints.
These local problems allow identifying the downscaling map
from average macroscopic quantities to the fine-grid variables.
By imposing the constraints for each continuum variable via the source term,
we define effective fluxes and the homogenized equation.
Using the local solutions in the oversampled regions with constraints
allows localizing the global downscaled map, which provides an accurate
representation of the solution; however, it is expensive as it involves
solving the global problem. By using the constraints in the oversampled
regions, we can guarantee the proximity between the global
and local downscaled maps for a given set of oversampled constraints.

The resulting homogenized
equation significantly differs from standard homogenization. First, there
are several variables per coarse block, which represent each continuum.
Secondly, the local problems are formulated in oversampled regions
with constraints.
Finally, the nonlinear homogenized fluxes depend on all averages
in oversampled regions, which bring non-local behavior for the equation.
These ingredients are needed to perform upscaling in the absence
of scale separation and high contrast.

In summary, our nonlinear nonlocal multicontinua approach has the following
steps.
\begin{itemize}

\item For each coarse-grid block, identify coarse-grid variables
(continua) and associated quantities. This is typically done via
some local spectral problems and describes the quantities
that can not be localized.

\item  Define local problems in oversampled regions with constraints.
The constraints are given for each continua variables over all coarse-grid
blocks. Local problems use specific source terms and boundary conditions,
which allow localizing the global downscaled map and
identifying effective fluxes.
This step gives a
downscaling from average continua to the fine-grid solution.

\item The formulation of the coarse-grid problem. We seek
the coarse-grid variables such that the downscaled fine-grid
solution approximately solves the global problem in a weak sense.
The weak sense is defined via specific test functions, which are
piecewise constants in each continua. This provides an upscaled model.

\end{itemize}

We consider two types of methods. First, we call linear interpolation
and the second is non-linear interpolation. In the first approach,
we seek the solution in the form
\[
w=\sum_{i,j} w_i^{j} \phi_i^{j},
\]
where $\phi_i^{j}$ are multiscale basis functions for coarse
cell $i$ and for continua $j$ defined in the oversampled region.
This approach is simpler as one uses linear approximation of the solution.
However, because of the nonlinearity, the nonlinear approximation is needed
(see \cite{ep03a,ep03c,ep03d}). In \cite{ep03d},
the authors propose such approach
for pseudomonotone equations (similar to homogenization).
In our problem, we present a nonlinear interpolation (we refer
as nonlinear nonlocal multicontinua). In this approach, the solution
is sought as a nonlinear map, which approximated $w_i^{j}$ from neighboring
cells in a nonlinear fashion, which is defined via local problems.
The local problem provides a nonlinear map from coarse-grid macroscopic
variables to the fine-grid solution.
This mapping is used to construct macroscopic equations.
We use non-locality
and multi-continua to address the cases without scale separation
and high contrast.

As one of the numerical examples, we consider two-phase flow and transport,
though our approach can be applied to a wide range of problems.
Two-phase flow and transport is one of challenging problems, where
many attempts are made to study it.
Typical approaches include pseudo-relative permeability approach,
also known as multi-phase upscaling, where the relative permeabilities
are computed based on local two-phase flow simulations. It is known
that these approaches are process dependent. Our approach provides a novel
upscaling, which shows that each coarse-grid block needs to contain several
average pressures and saturations and, moreover, these coarse-grid
relative permeabilities are non-local and depend on saturations and
pressures of neighboring cells. This model provides a new way to
look at two-phase flow equations and can be further used in
different modeling purposes.

The paper is organized as follows. In the next section, Section
\ref{sec:prelim}, we present some preliminaries and discuss
main concepts used for linear problems.
Section \ref{sec:concept} is devoted to a general concept
of nonlinear upscaling and this is discussed on several examples.
In Section \ref{sec:linear}, we present the linear interpolation
based approach and numerical results for multi-phase flow
and transport. Though the paper's main
focus is on nonlinear interpolation, the linear interpolation approach
is practical for many applications. In Section \ref{sec:nonlinear},
we present the nonlinear approach for multi-phase flow and
corresponding numerical results.

\section{Preliminaries}
\label{sec:prelim}

We consider a general nonlinear
system of the form
\begin{equation}
\label{eq:nonlin1}
\partial_t U + G(x,t,U,\nabla U) = g.
\end{equation}
In general, $U$ is a vector-valued function and
$g(x,t)$ is the source term. $G$ is assumed to be
heterogeneous in space (and time, in general),
which are resolved on the fine grid
(see Figure \ref{fig:ill1} for fine and coarse grid illustrations).
Our objective is to derive coarse-grid equations. First, we discuss
some examples for (\ref{eq:nonlin1}).

{\bf Example 1.} Nonlinear (pseudomonotone)
 parabolic equations, where
$G(x,t,U,\nabla U):=\text{div} \kappa(x,t,U,\nabla U)$.

{\bf Example 2.} Hyperbolic equations,
where $G(x,t,U,\nabla U): = v(x) \cdot \nabla G(x,U) $ and $U$ is
a scalar function and $v$ is a vector valued function.
For example, $v=-\kappa(x) \nabla p$ and $\text{div} (v ) =q$,
as in Darcy's flow.

{\bf Example 3}. Hamilton-Jacobi equations,
where $G(x,t,U,\nabla U): = G(x,\nabla U) $ and $U$ is
a scalar function.

{\bf Example 4.} One of our objective in this paper is to address the upscaling
of multi-phase flow, which is a challenging problem.
In this case, the model equations have the following form.
The two-phase flow equations are derived by writing Darcy's flow for each
water ($\alpha=w$) and oil ($\alpha=o$)
phases and the mass conservation as follows
\begin{equation}
\begin{split}
u_\alpha=- \lambda_{rw}(S_w) \kappa(x) \nabla p\\
\nabla \cdot (u_t) = q_p,\ \ \partial_t S_\alpha + \nabla \cdot (u_\alpha)=q.
\end{split}
\end{equation}
Here,
$u_{\alpha}$ is the Darcy velocity of the phase
$\alpha$, $u_t=\sum_\alpha u_\alpha$,
$S_{\alpha}$ is the saturation of the phase
$\alpha$, and
$\lambda_{r\alpha}$ is the relative mobility of the phase
$\alpha$.
By denoting the saturation of the water phase via $S=S_w$, we can
write the equations in the following way
\begin{equation}
\begin{split}
-div(\lambda(S) \kappa \nabla p) = q_p \\
\partial_t S + \nabla \cdot (u f(S))=q, \ \ u = - \kappa \nabla p.
\end{split}
\end{equation}
The special case of the two-phase flow consists of single-phase flow
\begin{equation}
\begin{split}
-div(\kappa \nabla p) = q_p \\
\partial_t S + \nabla \cdot (u S)=q, \ \ u = - \kappa \nabla p.
\end{split}
\end{equation}

\subsection{Overview of NLMC for linear problems}\label{sec:overview}

In this section, we will give a brief overview of the NLMC method
for linear problems \cite{NLMC}.
Our goal is to summarize the key ideas in linear NLMC and motivate the ideas in the extension of it to the nonlinear NLMC in the next section.
To be specific, we will consider
a model parabolic equation with a heterogeneous coefficient, namely,
\begin{equation}\label{eq:parabolic}
\frac{\partial u}{\partial t} - div(\kappa \nabla u) = g, \quad \text{ in } \Omega \times [0,T],
\end{equation}
with appropriate initial and boundary condition, where $\kappa$ is the heterogeneous field,
$g$ is a given source, $\Omega$ is the physical domain and $T>0$ is a fixed time.

\begin{figure}[!ht]
\centering
\includegraphics[scale=0.6]{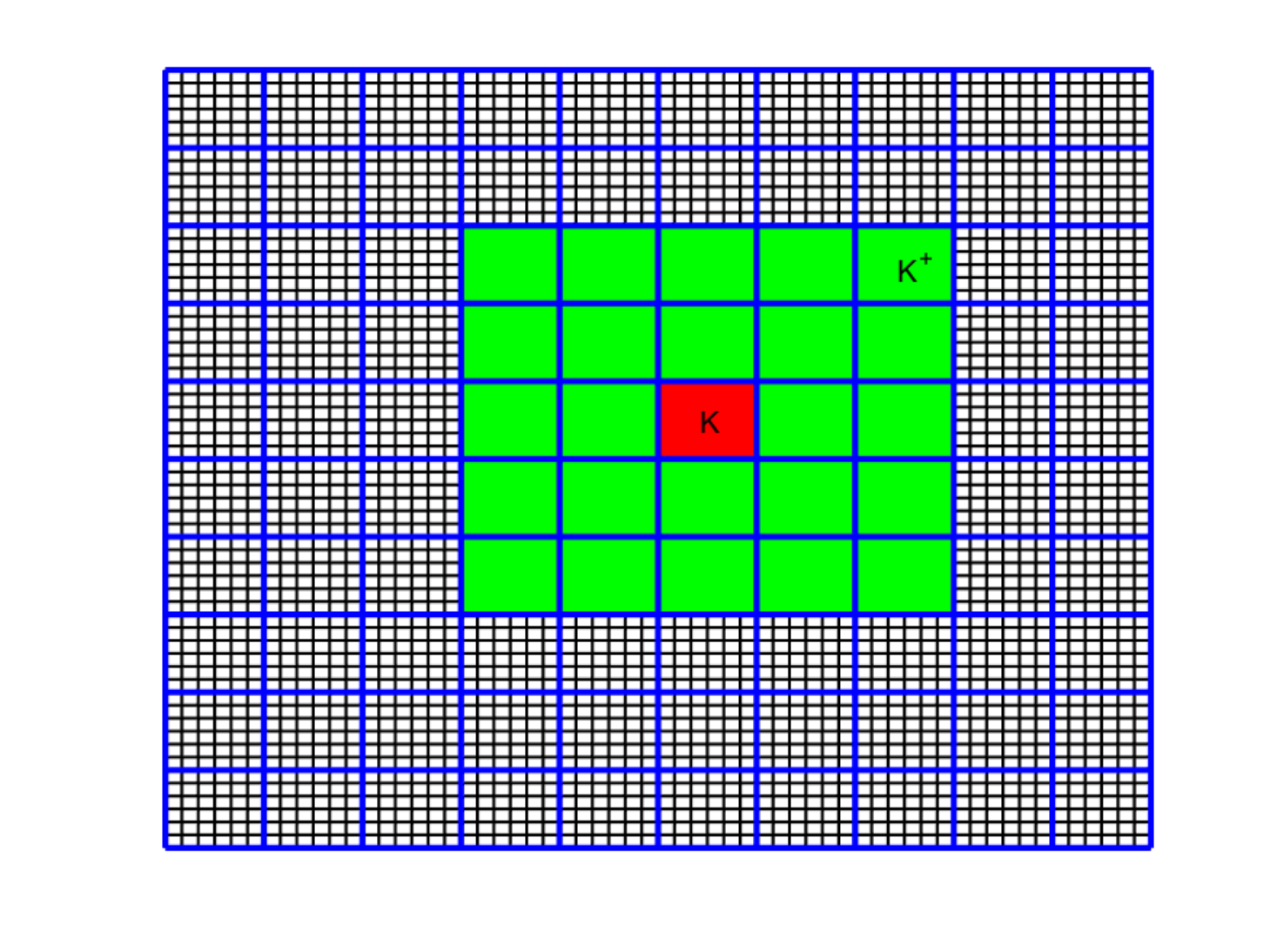}
\caption{Illustration of coarse and fine meshes as well as oversampled regions. The region $K^+$ (in Green)
is an oversampled region corresponding to the coarse block $K$ (in red).}
\label{fig:ill1}
\end{figure}

The NLMC upscaled system is defined on a coarse mesh, $\mathcal{T}^H$, of the domain $\Omega$.
We write $\mathcal{T}^H = \bigcup \{ K_i \; | \; i=1,\cdots, N\}$, where $K_i$ denotes the $i$-th coarse element
and $N$ denotes the number of coarse elements in $\mathcal{T}^H$.
For each coarse element $K_i$, we will identify multiple continua corresponding to various solution features. This can be done
via a local spectral problem or a suitable weight function.
The upscaled parameter is defined using multiscale basis functions.
For each coarse element $K_i$ and each continuum within $K_i$, we will construct
a multiscale basis function whose support is an oversampled region $K_i^+$, which is obtained by enlarging
the coarse block $K_i$ by a few coarse grid layers.
See Figure~\ref{fig:ill1} for an illustration of coarse grid and oversample region.
In particular, a structured coarse grid is shown with boundaries of coarse elements are denoted blue.
A coarse cell $K$ is denoted red and its oversampled region $K^+$ obtained by enlarging $K$ by two coarse grid layers
is denoted green.

Now we will specify the definition of continuum. For each coarse block $K_i$,
we will identify a set of continua which are represented by a set of auxiliary basis functions $\phi^{j}_i$,
where $j$ denotes the $j$-th continuum.
There are multiple ways to construct these functions $\phi^{j}_i$.
One way is to follow the idea proposed in CEM-GMsFEM \cite{chung2018constraint}.
In this framework, the auxiliary basis functions $\phi^{j}_i$
are obtained as the dominant eigenfunctions of a local spectral problem defined on $K_i$.
These eigenfunctions can capture the heterogeneities and the
contrast of the medium. The resulting solution of the upscaled system
corresponds to the moments of the true solution with respect to these
eigenfunctions.
Another way is to follow the framework in the original
NLMC method \cite{NLMC}, designed for flows in fractured media,
which can be easily modified for general heterogeneous media.
In this approach, one identifies explicit information of fracture networks.
The auxiliary basis functions $\phi^{j}_i$ are piecewise constant
functions, namely, they equal one within one fracture network
and zero otherwise.
The resulting solution of the upscaled system
corresponds to the average of the solution on fracture networks,
which can be regarded as one of the continua.
Finally, one can generalize the previous approach
to construct other auxiliary basis functions.
In particular, we can define $\phi^{j}_i$ to be characteristic functions
of some regions $K_{i}^{(j)} \subset K_i$.
These regions $K_{i}^{(j)}$ can be chosen to reflect various properties of the medium.
For instances, one can choose $K_{i}^{(j)}$ to be the region in which the medium has a certain range of values.
We will consider this last choice in this paper.



Once the auxiliary basis functions $\phi^{j}_i$ are specified,
we can construct the required basis functions.
The idea generalizes the original energy minimization framework in CEM-GMsFEM.
Consider a given coarse element $K_i$ and a given continuum $j$ within $K_i$.
We will use the corresponding auxiliary basis function $\phi^{j}_i$
to construct our required multiscale basis function $\psi^{j}_i$
by solving a problem in an oversampled region $K_i^+$.
Specifically, we find $\psi^{j}_i \in H^1_0(K_i^+)$ and $\mu \in V_{aux}$
such that
\begin{equation}\label{eq:basis}
\begin{aligned}
& \int_{K_i^+} \kappa \nabla \psi^{j}_i \cdot \nabla v + \int_{K_i^+} \mu v = 0, \quad \forall v\in H^1_0(K_i^+), \\
& \int_{K_\ell}\psi^{j}_{i}\phi_m^{\ell}   = \delta_{j\ell} \delta_{im}, \quad \forall K_\ell \subset K_i^+,
\end{aligned}
\end{equation}
where $\delta_{i m}$ denotes the standard delta function and $V_{aux}$ is the space spanned by auxiliary basis functions.
We remark the the function $\mu$ serves as a Lagrange multiplier for the constraints in the second equation of (\ref{eq:basis}).
We also remark that the basis function $\psi^{j}_i$ has mean value one on the $j$-th continuum within $K_i$
and has mean value zero in all other continua in all coarse elements within $K_i^+$.
In practice, the above system (\ref{eq:basis}) is solved in $K_i^+$ using a fine mesh, which is typically a refinement
of the coarse grid. See Figure~\ref{fig:ill1} for an illustration.


Now, we can derive the NLMC upscaled system.
As an illustration of the concept, we consider a forward Euler method for the time discretization of (\ref{eq:parabolic}).
At the $n$-th time step, the solution vector is denoted by $U^n$,
where each component of $U^n$ represents the average of the solution on a continuum within a coarse element.
We note that the size of this vector is $\sum_{i=1}^N L_i$, where $L_i$ is the number of continua in $K_i$.
The upscaled stiffness matrix $A_T$ is defined as
\begin{equation}\label{eq:Trans1}
(A_T)_{jm}^{(i,\ell)} = a(\psi_j^{i} ,\psi_m^{\ell} ) := \int_{\Omega} \kappa \nabla \psi^{i}_j \cdot \nabla \psi_m^{\ell},
\end{equation}
and the upscaled mass matrix $M_T$ is defined as
\begin{equation}\label{eq:Trans2}
(M_T)_{jm}^{(i,\ell)}  = \int_{\Omega}  \psi^{i}_j  \psi_m^{\ell}.
\end{equation}
Finally, the NLMC system is written as
\begin{equation}
\label{eq:nlmc-linear}
M_T (U^{n+1} - U^n) + \Delta t A_T U^n = \Delta t G^n
\end{equation}
where the components of the vector $G^n$ are defined as $\int_{\Omega} g(\cdot, t_n) \phi^{j}_i$
and $t_n$ is the time at the $n$-th time step.
We remark that the nonlocal connections of the continua are coupled by the matrices $A_T$ and $M_T$.
We also remark that the local computation in (\ref{eq:basis})
results from a spatial decay property of the multiscale basis function, see \cite{chung2018constraint,chung2018fast,chung2018constraintmixed}
for the theoretical foundation.

\section{Nonlinear non-local multicontinua model}

\label{sec:concept}
\subsection{General concept}

We will first present some general concept of our nonlinear NLMC.
We consider the following nonlinear problem
\begin{equation}
\label{eq:non1}
U_t + G(x, U, \nabla U)=g,
\end{equation}
where
$G$ has a multiscale dependence with respect to
space (and time, in general).

To formulate the coarse-grid equations in the time interval
$[t_n,t_{n+1}]$, we first introduce coarse-grid variables
$U_{i}^{n,j}$, where $i$ is the coarse-grid block, $j$ is
a continuum representing the coarse-grid variables,
and $n$ is the time step (cf. \cite{NLMC}).
As we noted that the continuum plays a role of a macroscopic variable,
which can not be localized.
For each coarse-grid block $i$, we need several
coarse-grid variables, which will be indexed by $j$.
In this paper, we will consider two types of interpolation
for multiscale degrees of freedom.
The first approach is called linear interpolation, which constructs approximate solution
$U_H^n$ at the time $t_n$ as
\[
U_H^n = \sum_{i,j} U_i^{n,j} \psi_i^{j},
\]
where $\psi_i^{j}$ are multiscale basis functions, which are defined
via local constrained problems formulated in the oversampled
regions. These basis functions are possibly
supported in the oversampled regions
and the resulting coarse-grid equations will be nonlinear and non-local
when projected to the appropriate test spaces. The
values of $U_i^{n,j}$ are computed via the variational formulation
using an explicit or implicit discretization of time
\begin{equation}
(U_H^{n+1}, V_H) - (U_H^{n}, V_H)
+ \Delta t (G(x,U_H^{L},\nabla U_H^{L}), V_H) = \Delta t (g,V_H),\ \ \forall V_H,
\end{equation}
where $L=n$ or $L=n+1$ depending whether explicit or implicit
discretization is used, and $V_H$ denotes test functions.
Here, $(\cdot,\cdot)$ denotes usual $L^2$ inner product.
 As for the test space, one can use
the multiscale space spanned by $\psi_i^{j}$ or the auxiliary basis functions $\phi_i^{j}$.
This linear approach is simpler to use; however, it ignores the nonlinear
interpolation, which is important for nonlinear homogenization and
numerical homogenization \cite{pankov97, ep03d}.

Next, we describe the nonlinear approach, which we refer
as nonlinear nonlocal multicontinuum approach. In this case,
the solution is sought as a nonlinear interpolation of the degrees
of freedom $U_i^{n,j}$ defined in the oversampled region.
More precisely,
we compute a downscale function using the given solution values $U_i^{n,j}$
by solving a constrained problem in the oversampled regions.
Mathematically, this downscale function $U^n_h$ can be written as a function of all values $U_i^{n,j}$, namely,
\[
U_h^n=\mathcal{F}(U^n),
\]
where $\mathcal{F}$ is a nonlinear map and $U^n$ is a vector containing all values $\{ U_i^{n,j}\}$.
The nonlinear map is computed by solving the local
problem (cf. (\ref{eq:local1})),
and the function $U^n_h$ is glued together from local downscaled maps.
Once this map is defined, we seek
all coarse-grid macroscopic variables $\{ U_i^{n,j}\}$ such that the downscaled
fine-grid solution $U^n_h$ solves the global problem in a variational setting.
The test functions are defined
for each degrees of freedom in the
form of piecewise constants or piecewise
functions, in for example a finite volume or Petrov Galerkin setting.
More precisely,
the coarse-grid system has the form
\begin{equation}
(U_h^{n+1}, V_H) - (U_h^{n}, V_H)
+ \Delta t
 (G(x,U_h^L,\nabla U_h^L), V_H) = \Delta t (g,V_H),
\end{equation}
for all suitable test functions $V_H$.
Here, $L=n$ or $L=n+1$  can be used for explicit or
implicit discretization. The test functions are chosen
to be piecewise polynomials or auxiliary basis functions $\phi^j_i$.
We remark that the dimension of the test space is chosen to be the dimension of the
coarse-grid macroscopic variables.

\subsection{Nonlinear nonlocal multicontinuum approach}
\label{sec:non-nlmc}

\begin{figure}[!ht]
\centering
\includegraphics[scale=0.6]{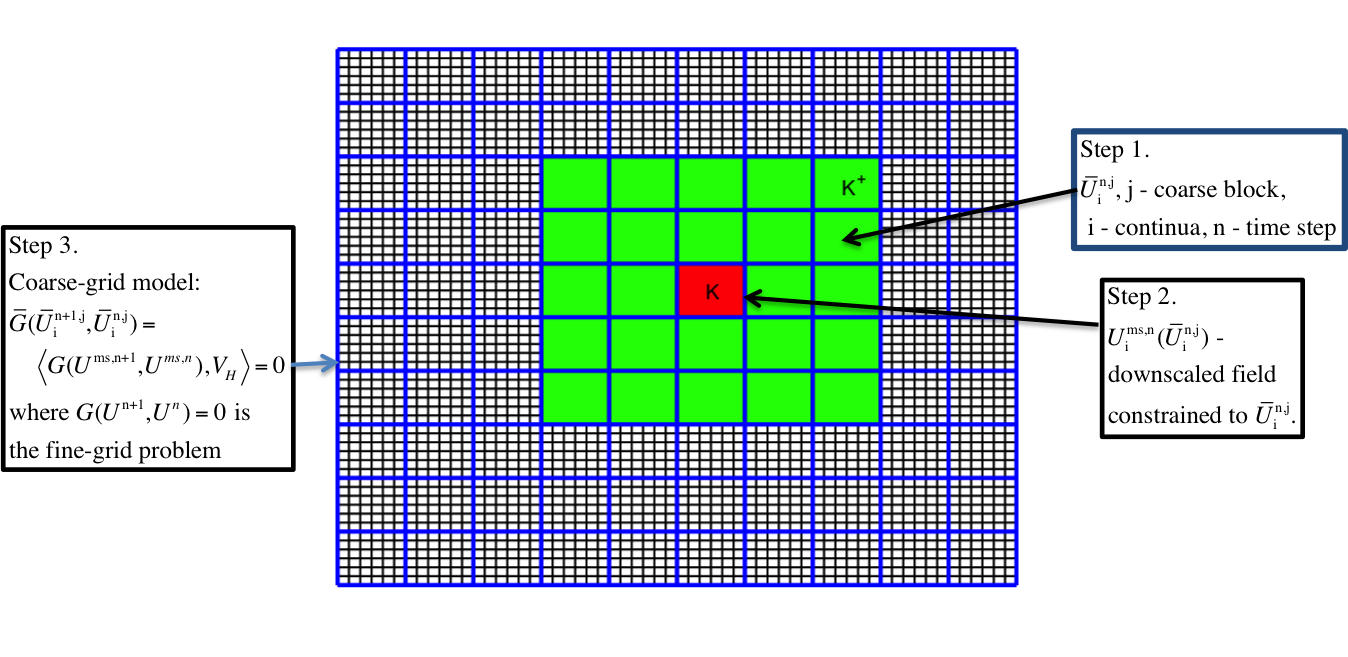}
\caption{Illustration of the steps.}
\label{fig:ill2}
\end{figure}

In this section, we give some details of our nonlinear nonlocal
multicontinuum approach, in general, and then present some
examples. The illustration is given
in Figure \ref{fig:ill2}.
In later sections, we present a detailed application to
two-phase flow equations.

{\bf Step 1. Defining coarse-grid variables.}
The first step involves defining macroscopic variables.
In our cases, we will define them by $U_{j}^{n,i}$,
where $i$ is the coarse-grid block, $j$ is the continua,
and $n$ is the time step.  These variables in our applications
will be defined by prescribing averages to subregions, which can have
a complex shapes (for example, fracture regions).
The equation for $U_{j}^{n,i}$, in general, will have a form
\begin{equation}
\label{eq:non1-discrete}
\overline{U}_j^{n+1,j} - \overline{U}_j^{n,j}  +
\overline{G}_{i}^{j}(\overline{U}^{L})=0,
\end{equation}
where $L=n$ or $L=n+1$
for explicit or implicit discretizations. The operator $\overline{G}_{i}^{j}$ is determined respect to
the given continuum $j$ and coarse element $i$, and contains information from the source term and the time step size $\Delta t$.
In (\ref{eq:non1-discrete}), we use $\overline{U}^L$ to denote the vector consisting of all macroscopic quantities $U_i^{L,j}$ for all $i$ and $j$.
The operator $\overline{G}_{i}^{j}$ is obtained by solving local problems on an oversampled region corresponding to $K_i$.
To define (\ref{eq:non1-discrete}), we next discuss local problems.


%
%
%

{\bf Step 2. Local solves for generic constraints.}
The computation of $\overline{G}^j_i$ requires solutions of constrained
local nonlinear problems.
Here we present a general formulation. We consider $i$ to be the index for a coarse element $K_i$
or a coarse neighborhood $\omega_i$, which is defined for a coarse node $i$ by
\[
\omega_i = \bigcup \{K\in \mathcal{T}^H:x_i\in \overline{K}\}.
\]
The choice of $K_i$ or $\omega_i$ depends on the global coarse grid discretization.
For example, with finite volume or Petrov Galerkin formulation, we choose $K_i$,
while for continuous Galerkin formulation, we choose $\omega_i$.
In this part, we present the steps using $\omega_i$.
The required local nonlinear problem will be solved on $\omega_i^+$,
which is an oversampled region obtained by enlarging $\omega_i^+$ a few coarse grid layers.
We let $c:=\{ c_m^{(l)} \}$ be a set of scalar values, where $m$ denotes the $m$-th coarse element in the oversampled region $\omega_i^+$
and $l$ denotes the $l$-th continuum within $K_m \subset \omega_i^+$.
The local problem is solved by finding a function $N_{\omega_{i}}(x;c)$ given by
\begin{equation}
\label{eq:local1}
G(x,N_{\omega_{i}}(x;c),\nabla N_{\omega_{i}}(x;c))  = \sum_{m,l} \mu_{i.m}^{(l)}(c) I_{K_{m}^{(l)}} \;\text{in }\omega_{i}^{+}
\end{equation}
with constraints
\[
\int_{\omega^+} N_{\omega_{i}} (x;c) I_{K_m^{(l)}}(x) = c_m^{(l)}.
\]
Here, $I_{K_m^{(l)}}(x)$ is an indicator function for the region
$K_m^{(l)}$ defined for the continua $l$ within coarse block $m$.
We remark that the values $\mu_{i.m}^{(l)}(c)$ play the role of Lagrange multipliers for the constraints.
We notice that $N_{\omega_i}$ depends  on all constraints in the oversampled region.
In general, one can also impose constraints on the gradients
of $N_{\omega_i}$ and the constraint function can have a complex form.
We remark that the problem (\ref{eq:local1}) requires a boundary condition, which is problem dependent. We will discuss it later.
Our main message is that the local problems with constraints are
solved in oversampled region and involves several
coarse-grid variables for each coarse block.
For the precise formulation, the values $\{ c_m^{(l)} \}$ will be chosen as the macroscopic variables.

{\bf Discussion on localization.}

\begin{figure}[!ht]
\centering
\includegraphics[scale=0.6]{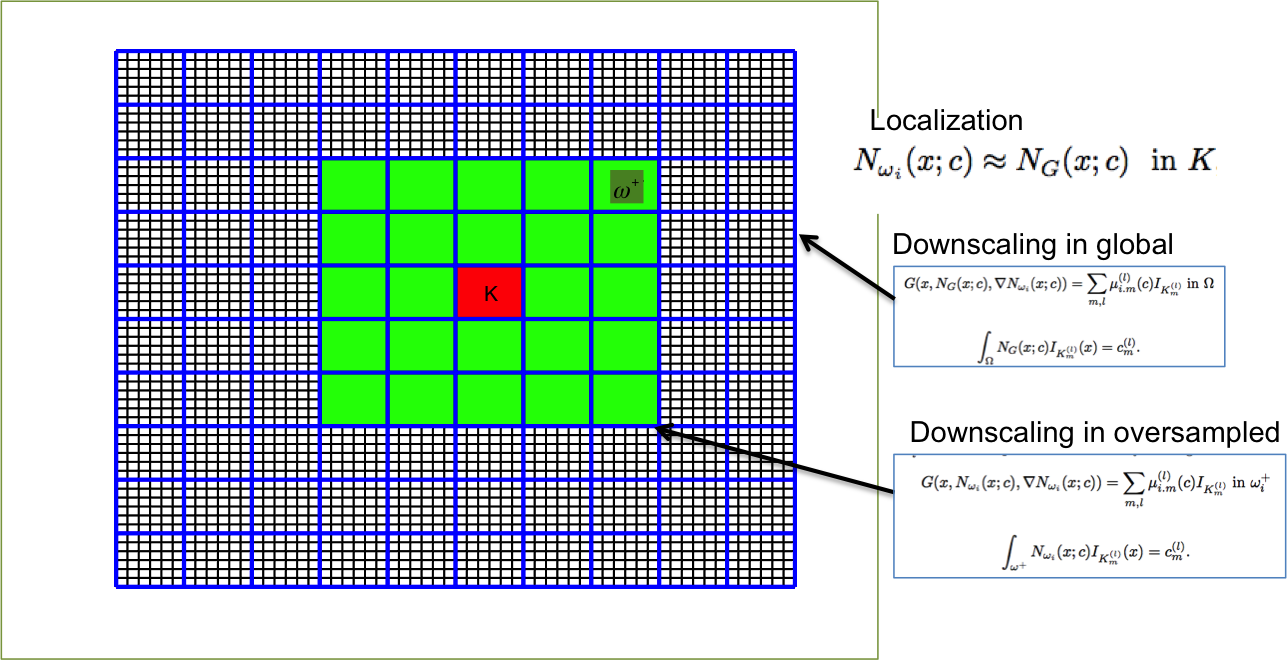}
\caption{Illustration of the localization.}
\label{fig:ill11}
\end{figure}

Next, we discuss the idea behind the localization principle stated
above (see Figure \ref{fig:ill11} for illustration). For this reason, we first
 define the global downscaling operator,
$N_{G}(x;c)$, which solves the global problem with
constraints
\begin{equation}
\label{eq:local1_G}
G(x,N_{G}(x;c),\nabla N_{G}(x;c))  = \sum_{m,l} \mu_{i.m}^{(l)}(c) I_{K_{m}^{(l)}} \;\text{in }\Omega
\end{equation}
with constraints
\[
\int_{\Omega} N_{G} (x;c) I_{K_m^{(l)}}(x) = c_m^{(l)}.
\]
The boundary conditions are taken as the original global boundary conditions.
For our localization, we desire that the local downscaled solution,
$N_{\omega_{i}}(x;c)$ in $K$ (target coarse block)
defined in (\ref{eq:local1}),
is approximately the same as $N_{G}$
restricted to $K$ for the same constraint values as the global
problem in $\omega_{i}^{+}$. I.e.,
\[
N_{\omega_{i}}(x;c)\approx N_{G}(x;c)\ \ \text{in } K,
\]
where $c$ values coincide in $\omega_{i}^{+}$.
This property allows solving local nonlinear problems instead
of the global problem without
sacrificing the accuracy
and defining correct upscaled coefficients.
This desired property can be shown for monotone elliptic
equations using zero Dirichlet boundary conditions on
$\partial \omega_i^{+}$.
For more complex nonlinear system, one needs to also use
appropriate boundary conditions on
$\partial \omega_i^{+}$ to achieve this property.

{\bf Step 3. Defining coarse-grid model.}
The upscaled
flux $\overline{G}_{i}^{j}$ is defined by substituting the downscaled
global solution corresponding to the constraints given by $c$ and
multiplying it by test functions.
For this procedure, one first  defines a global fine-grid
downscaled field that is constrained. The equation
(\ref{eq:local1}) defines local fine-grid fields constrained to
macroscopic variables in oversampled regions. In general,
one needs to ``glue'' together a global downscaled solution, which
approximates the global problem in a weak sense.
A simple approach to glue together is to use partition of unity functions,
$\chi_{\omega_i}$ for each region $\omega_i$. Then,
\[
\mathcal{F}(\overline{U})=\sum_i N_{\omega_i} \chi_{\omega_i}.
\]
In some discretization
 scenarios (e.g., finite volume, Discontinuous Galerkin,...),
it is not necessary to glue together
in order to
 obtain a global fine-grid field.
In the proposed NLMC approach, our coarse-grid fomrulation
uses ``finite volume'' type approximation and we will not need to glue
local approximations. In this case,
\begin{equation}
\label{eq:FV_approx}
\mathcal{F}(\overline{U})= N_{\omega_i} \text{ and } \nabla \mathcal{F}(\overline{U})= \nabla N_{\omega_i} \text{ on } E_i .
\end{equation}
where $E_i$ is a coarse edge for the coarse grid.
In general,
we can write it as
$\mathcal{F}(\overline{U})= N_{\omega_i}$
in all local fine-grid blocks, and the following variational formulation
\begin{equation}
\label{eq:nlnlmc}
({\partial \over\partial t}\mathcal{F}(\overline{U}),V_H) +(G(x,\mathcal{F}(\overline{U}),\nabla \mathcal{F}(\overline{U})),V_H) =  (g,V_H),
\end{equation}
where $V_H$ are test functions. In this setup,
\[
(G(x,\mathcal{F}(\overline{U}),\nabla \mathcal{F}(\overline{U})),V_H) =(\overline{G}(x,\overline{U}),V_H).
\]
The above equation can formally
be thought
as
\[
{\partial \over\partial t} \mathcal{M}(\overline{U}) + \mathcal{G} (\overline{x},\overline{U})=\overline{g}.
\]
The time discretization of Equation (\ref{eq:nlnlmc}) can be written as
\begin{equation}
\label{eq:nlnlmc1}
(\mathcal{F}(\overline{U}^{n+1}),V_H)-(\mathcal{F}(\overline{U}^{n}),V_H) +
\Delta t(G(x,\mathcal{F}(\overline{U}^{L}),\nabla \mathcal{F}(\overline{U}^L)),V_H)
 =  \Delta t (g,V_H),
\end{equation}
where $L=n$ or $L=n+1$ for explicit or implicit discretizations.
In general, the downscaling operator has
the following property
\begin{equation}
\label{eq:aver1}
(\mathcal{F}(\overline{U}^{L}),V_H) \approx
(\overline{U}^{L},V_H),
\quad\text{or}\quad
(\mathcal{F}(\overline{U}^{L}),V_H) =
(\overline{U}^{L},V_H)
\end{equation}
which simplifies the computations. In the latter situation, the mass matrix $\mathcal{M}$ is diagonal.
In addition,
\begin{equation}
\label{eq:aver11}
(\mathcal{F}(\overline{U}^{L}),V_h)\approx
({U}_h^{L},V_h),
\end{equation}
where ${U}_h^{L}$ is a fine-scale field and
$\overline{U}^{L}$ is the corresponding coarse-grid average.
The equation (\ref{eq:aver11}) states that if we use the average
of the fine-scale field to re-construct it, the resulting approximation
remains close to the original fine-scale field.
We remark that the operator $\overline{G}_i^j$ in (\ref{eq:non1-discrete})
is obtained using (\ref{eq:nlnlmc1}). More precisely, we have
$$\overline{G}_i^j :=
\Delta t(G(x,\mathcal{F}(\overline{U}^{L}),\nabla \mathcal{F}(\overline{U}^L)),V_H)
 -  \Delta t (g,V_H)$$
 by using a test function $V_H$ corresponding to region $i$ and continuum $j$.

Next, we present some more concrete constructions for some cases, listed in Section~\ref{sec:prelim}.

{\bf Example 1.} In this example,
we consider a class of nonlinear (pseudomonotone)
 parabolic equations with
$G=\text{div} \Big( k(x,U,\nabla U)\Big)$, where $k$ is a heterogeneous function in $x$, and in general, also in $U$ and $\nabla U$.
We will consider a continuous Galerkin discretization in space on a coarse grid and explicit Euler discretization in time.
Let $c = \{ U^{n,j}_i\}$ be a set of macroscopic variables at the time $t_n$. In particular, $U^{n,j}_i$ denotes the mean value of the solution on continuum $j$ within coarse element $i$
at time time $t_n$.
On an oversampled region $\omega_i^+$, we find a function $N_{\omega_i}$ on the fine grid such that
\[
G(x,N_{\omega_i},\nabla N_{\omega_i}) = \sum_{K_m^{(l)} \subset \omega_i^+} \mu_{m}^{(l)}(c)I_{K_{m}^{(l)}} \;\quad\text{in } \;\omega_i^+
\]
subject to the following
constraints
\begin{equation*}
\begin{split}
&\int_{\omega_i}  N_{\omega_i}(x) I_{K_m^{(l)}}(x) = U^{n,l}_m, \quad \forall K_m^{(l)} \subset \omega_i, \\
&\int_{\omega_i}  N_{\omega_i}(x) I_{K_m^{(l)}}(x) = 0, \quad \forall K_m^{(l)} \subset \omega_i^+\backslash\omega_i.
\end{split}
\end{equation*}
We remark that $\mu_{m}^{(l)}$ plays the role of Lagrange multiplier. The above problem is solved using the Dirichlet boundary condition on $\partial\omega_i^+$.
We can set $N_{\omega_i}$ equals to zero on $\partial\omega_i^+$, and this choice of motivated by the decay property of multiscale basis functions in CEM-GMsFEM.
Then we can define a global fine scale function by
$$
U^n_h = \sum_i N_{\omega_i} \chi_{\omega_i}
$$
where $\{ \chi_{\omega_i}\}$ are suitable partition of unity functions, where we notice that the values of $N_{\omega_i}$ within $\omega_i^+\backslash\omega_i$
are not used.
Using the global downscale function $U^n_h$ and suitable test functions,
we can derive a coarse grid scheme.
We note that the upscaled model has the form (\ref{eq:nlnlmc}).

{\bf Example 2.}
In this example, we consider a class of
hyperbolic equations,
where $G(x,t,U,\nabla U): = v(x) \cdot \nabla G(x,U) $ and $U$ is
a scalar function and $v$ is a vector valued function.
For example, $v=-\kappa(x) \nabla p$ and $\text{div} (v ) =q$,
as in Darcy's flow.
We will consider a finite volume type discretization in space on a coarse grid and explicit Euler discretization in time.
Let $c = \{ U^{n,j}_i\}$ be a set of macroscopic variables at the time $t_n$. Same as above, $U^{n,j}_i$ denotes the mean value of the solution on continuum $j$ within coarse element $i$
at time time $t_n$.
On an oversampled region $K_i^+$, we find a function $N_{K_i}$ on the fine grid such that
\[
G(x,\nabla N_{K_i}) = \sum_{K_m^{(l)} \subset K_i^+} \mu_{m}^{(l)}(c)I_{K_{m}^{(l)}} \;\quad\text{in } \; K_i^+
\]
subject to the following
constraints
\[
\int_{K_i^+}  N_{\omega_i}(x) I_{K_m^{(l)}}(x) = U^{n,l}_m, \quad \forall K_m^{(l)} \subset K_i^+.
\]
We remark that $\mu_{m}^{(l)}$ plays the role of Lagrange multiplier. The above problem is solved using the standard inflow boundary condition on the inflow part of $\partial K_i^+$.
We can set $N_{K_i}$ equals to zero to the value $\{ U^{n,j}_i\}$ chosen by upwinding.
Then we can define a global fine scale function by
$$
U^n_h = \sum_i N_{K_i} \chi_{K_i^+}
$$
where $\{ \chi_{K_i^+}\}$ are suitable partition of unity functions for the overlapping partition $\{ K_i^+\}$.
We remark that this global downscale function $U^n_h$ preserves the mean values, namely,
$$
\frac{1}{|K_i^{(j)}|}\int_{K_i^{(j)}} U^n_h = U^{n,j}_i.
$$
Using the global downscale function $U^n_h$ and suitable test functions,
we can derive a coarse grid scheme.
For the test functions, we use $\{ I_{K_m^{(l)}}(x)\}$, for all $m,l$. Hence, we obtain
\begin{equation}
\label{eq:up}
U^{n+1,j}_i = U^{n,j}_i - \Delta t \int_{K_i^{(j)}} G(x,\nabla U^n_h) + \Delta t \int_{K_i^{(j)}} g, \quad\quad \forall \, i,j.
\end{equation}
We note that this upscaled model has the form (\ref{eq:nlnlmc}).
We also note that the number of unknowns is the same of number of equations.

{\bf Example 3.}
For this case, the local problem can be formulated as follows.
On an oversampled region $K_i^+$, we find a function $N_{K_i}$ on the fine grid such that
\[
G(x,\nabla N_{K_i}) = \sum_{K_m^{(l)} \subset K_i^+} \mu_{m}^{(l)}(c)I_{K_{m}^{(l)}} \;\quad\text{in } \; K_i^+
\]
subject to the following
constraints
\[
\int_{K_i^+}  N_{\omega_i}(x) I_{K_m^{(l)}}(x) = U^{n,l}_m, \quad \forall K_m^{(l)} \subset K_i^+.
\]
Similar to Example 2, the upscaled model has the form (\ref{eq:nlnlmc}) and (\ref{eq:up}).

In summary, the local problems involve original local problems
with constraints formulated in the oversampled regions.
The source term constants represent the homogenized fluxes and their
dependence on averages of solutions and gradients are the functional form
of the equations.

\subsubsection{RVE-based extension of the method}

The proposed concept can also be used for problems with scale separation.
A typical example is a fractured media, where the fracture distributions
are periodic or posess some scale separation.
In this case, we do not use the oversampling based local problems
and simply use local problems in RVE. To be more precise,
Step 1 remains the same as before, which involves identifying
local multicontinua variables in each coarse-grid block.

In Step 2, we solve the local problem (\ref{eq:local1})
in RVE subject to the constraints. The main difference is that
we use only the constrained in the target coarse-grid block. Thus,
there are fewer constraints and the problem is localized to
the target coarse block.

Once the local solves are defined, via periodicity, we extend the local
solution to the $\omega_i$ and use this extension to compute
the effective flux. This construction is similar to e.g.,
\cite{ep03d} or MsFEM using RVE \cite{eh09}.

\section{Linear Approach}
 \label{sec:linear}

 In this section, we will focus on the linear approach (c.f. Section~\ref{sec:concept}),
 and present some numerical results.

\subsection{Linear transport}
\label{sec:lin-tran}

\begin{figure}[!ht]
\centering
\includegraphics[scale=0.4]{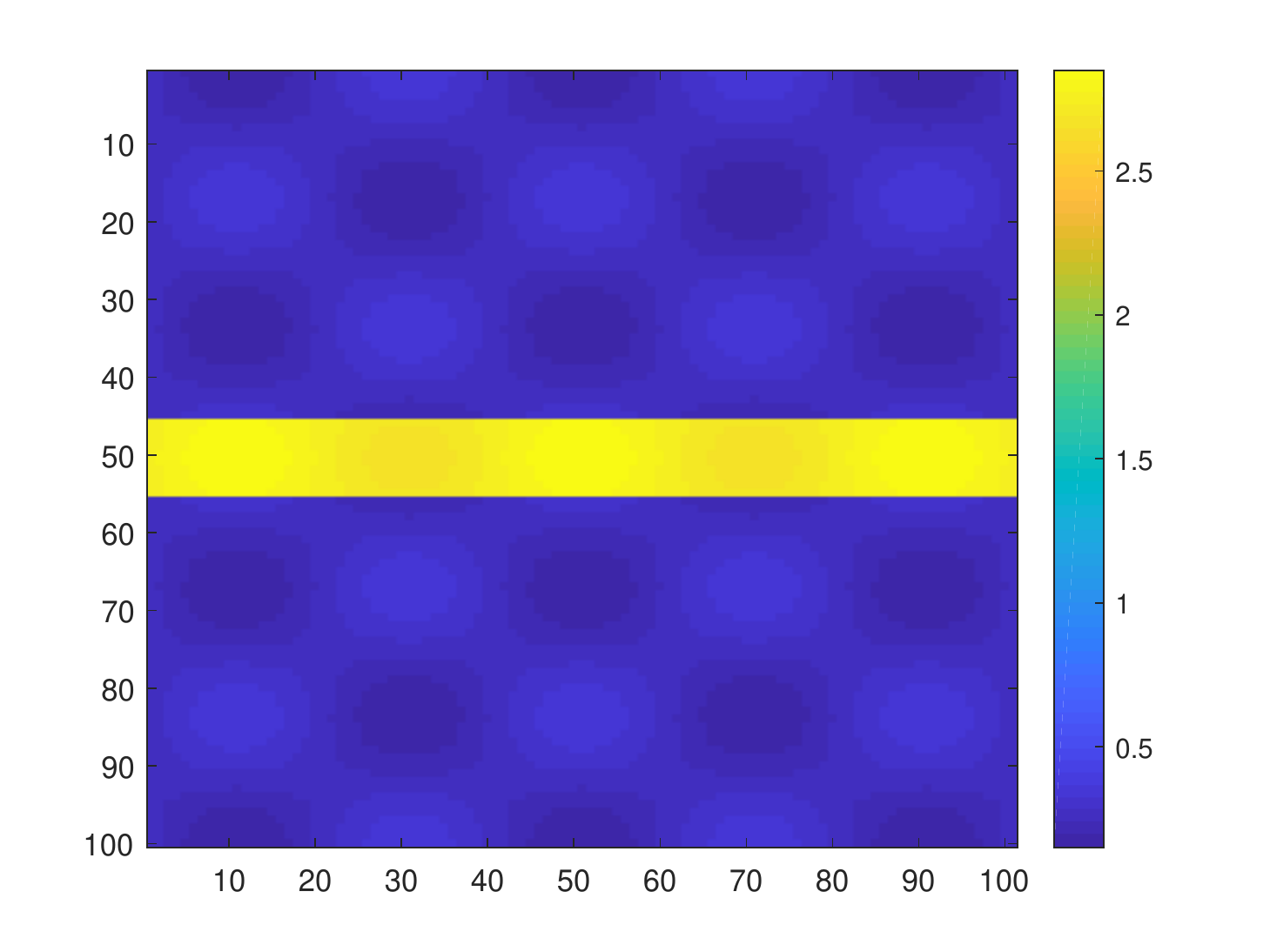} \includegraphics[scale=0.4]{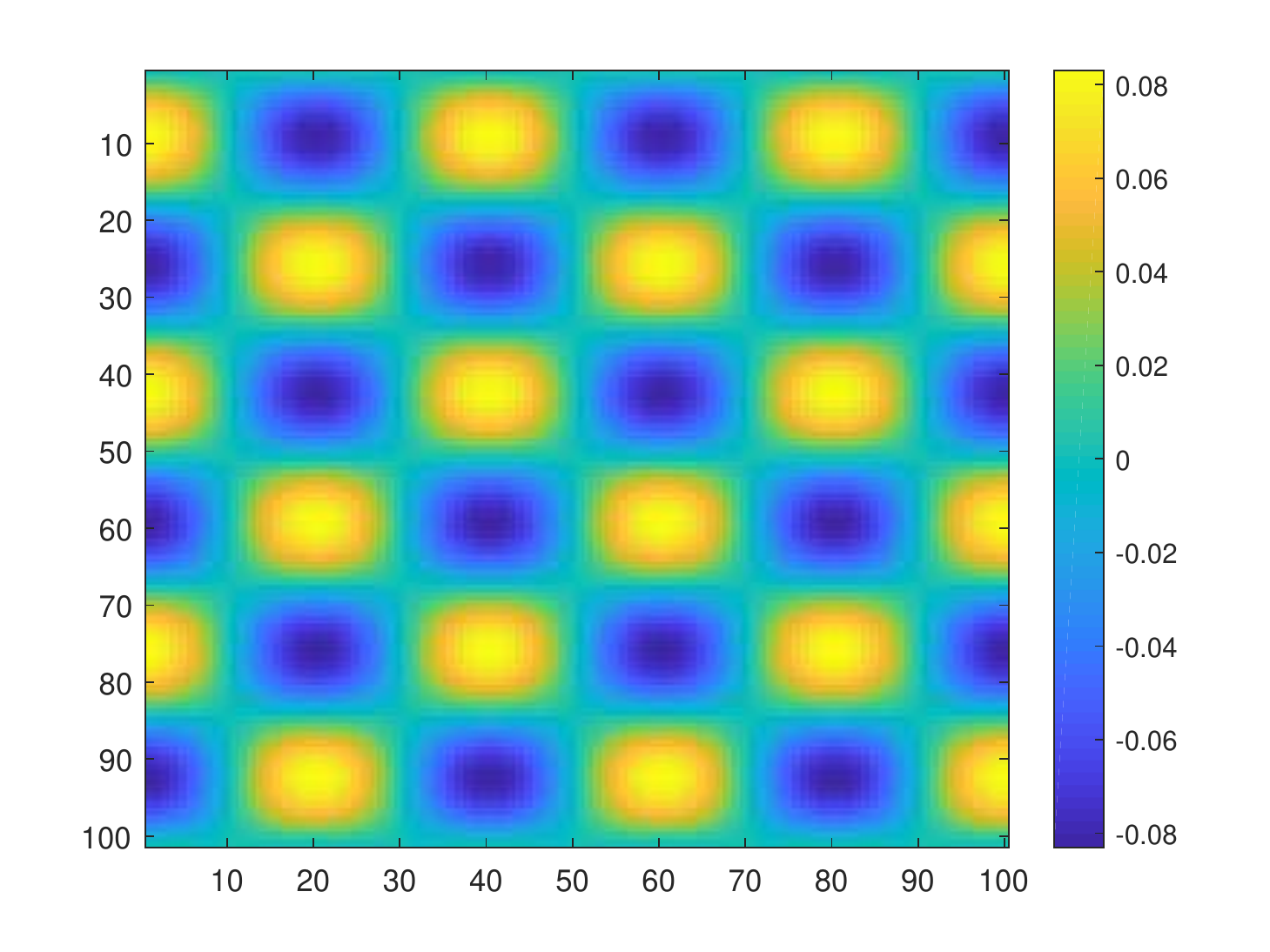}
\caption{A velocity field $u$. Left: $x_{1}$-component. Right: $x_{2}$-component.}
\label{fig:velocity_case1}
\end{figure}

In this section, we will present some numerical results
for the linear transport equation with a given velocity.
More precisely, we consider
\[
\partial_{t}S+\nabla\cdot(uS)=q\;\text{in }\Omega
\]
where $u\in H(\text{div},\Omega)$ is a given divergence-free velocity field with
$u\cdot n=0$ on $\partial\Omega$.
We take $\Omega=[0,1]^{2}$, and
the velocity field $u$ is shown in Figure \ref{fig:velocity_case1}.
The coarse grid size $H=1/20$, and we use a single continuum model.
In Table \ref{tab:trans_error_case1}, we present a error comparison
between our upscaling method and the standard finite volume method.

The derivation of the upscaled system follows the ideas in Section~\ref{sec:overview}. For each coarse element $K_i$,
we solve the following system in an oversampled region $K_i^+$:
\begin{equation}
\label{eq:lin-trans}
\nabla \cdot (u\psi_i) + \mu = 0
\end{equation}
to obtain a basis function $\psi_i$. The above system is equipped with the constraints
\begin{equation}
\int_{K_j} \psi_i = \delta_{ij}.
\end{equation}
We remark that $\mu$ is a piecewise constant function, and plays the role of Lagrange multiplier.
We notice that the upscaled system has the form (\ref{eq:nlmc-linear}).
In practice, one can add artificial diffusion in (\ref{eq:lin-trans}) in order to obtain a stable numerical solution,
and use zero Dirichlet boundary condition. We remark that one can use other boundary conditions, see Section~\ref{sec:non-nlmc}.

\begin{figure}[!ht]
\centering
\includegraphics[scale=0.35]{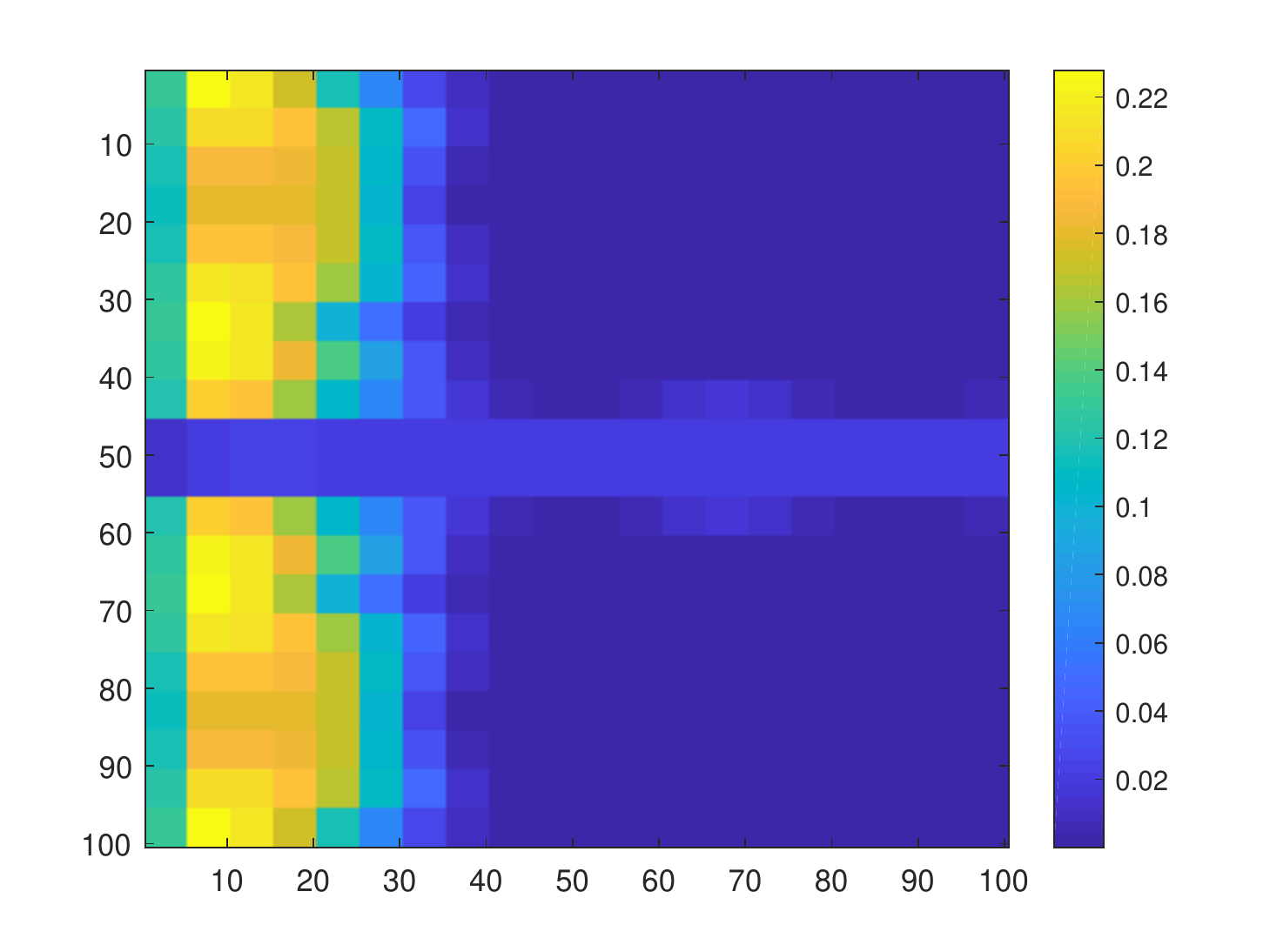}
\includegraphics[scale=0.35]{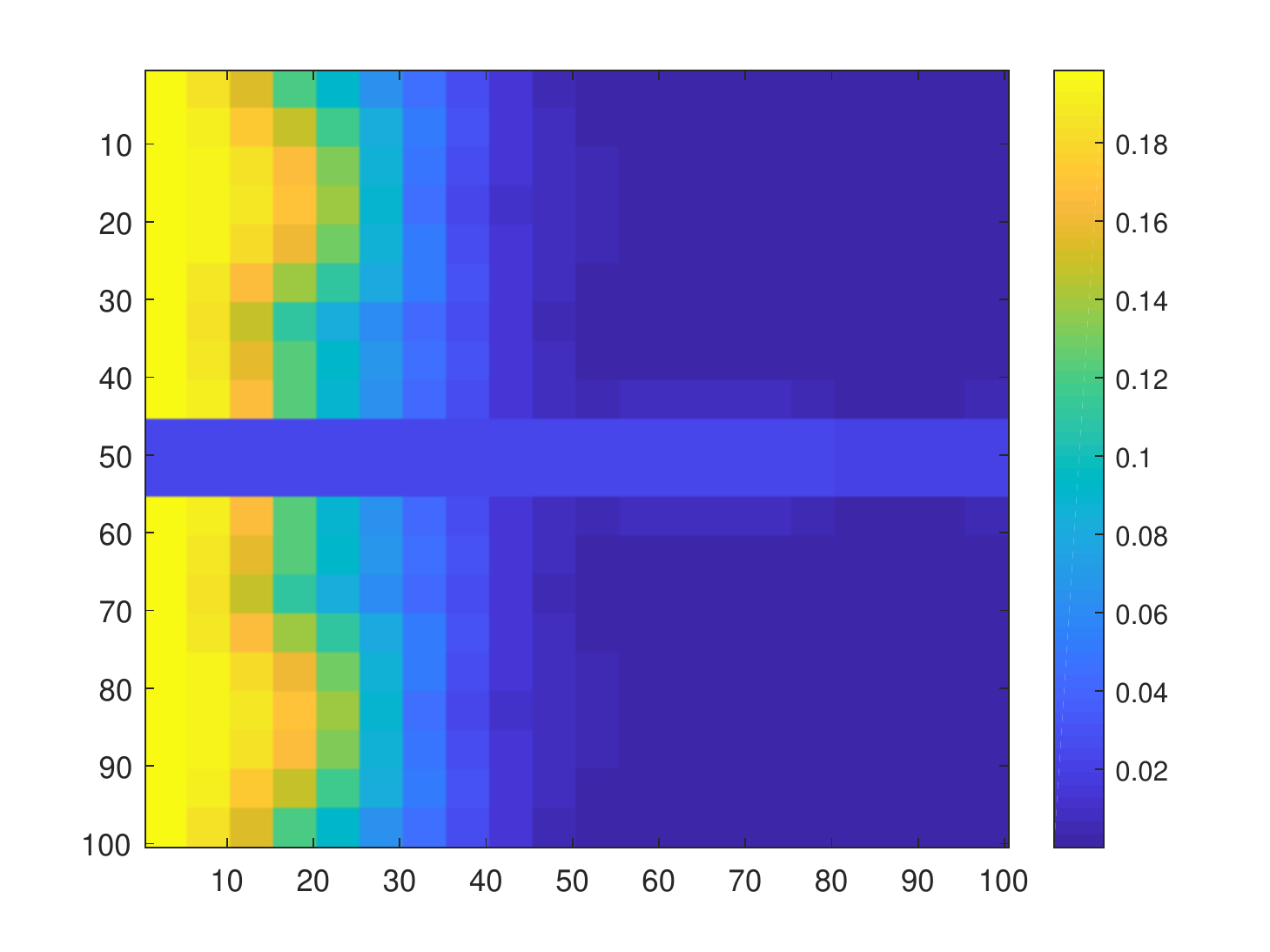}
\includegraphics[scale=0.35]{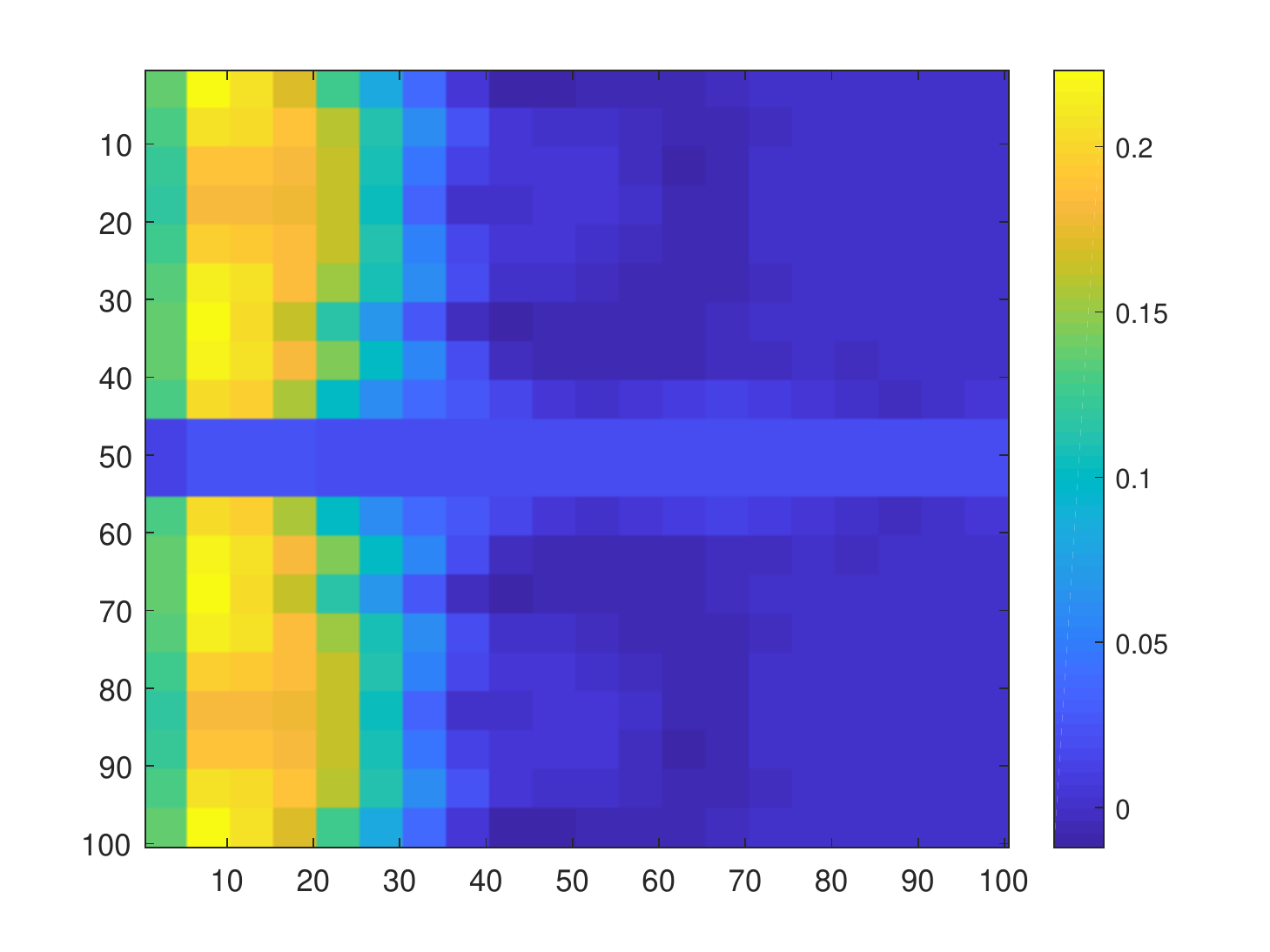}
\caption{Snapshots of the solution at $T=1$. Left: fine solution, Middle: finite volume solution,
Right: upscale solution with $5$ oversampling layers.}
\label{fig:snap1}
\end{figure}

By using the standard finite volume method, the relative $L^2$ error is $27.78\%$ at the final time $T=1$.
From Table \ref{tab:trans_error_case1}, we see that our upscaling method provides much better approximations,
where "$\#$ layer" standards for the number of layers used in the oversampling domains.
When $\#$ layer equals $\infty$, it means that the oversampling domain is the whole physical domain.
In Figure~\ref{fig:snap1}, we show
the snapshots of the solution at $T=1$.  From this figure, we observe that
the solution computed with CEM provides a good approximation, particularly,
in saturated regions.


\begin{table}[!ht]
\centering
\begin{tabular}{|c|c|}
\hline
\#layer  & Errors \tabularnewline
\hline
\hline
5 & 11.57\%\tabularnewline
\hline
7 & 7.73\%\tabularnewline
\hline
9 & 7.22\%\tabularnewline
\hline
$\infty$ & 6.91\%\tabularnewline
\hline
\end{tabular}
\caption{$L^2$ relative errors for our upscaling method. The relative error for standard finite volume scheme on the same grid is $27.78\%$.}
\label{tab:trans_error_case1}
\end{table}

\subsection{Single phase flow}

In this section, we apply our method to a single phase flow problem
which is described as follows:
\begin{align}
\partial_{t}S_{w}+\nabla\cdot(u_{w}S_{w}) & =q_{w}, \label{eq:single_phase_eq1}\\
-\nabla\cdot(\kappa\nabla P_{w}) & =q, \label{eq:single_phase_eq2}\\
u_{w} & =-\kappa\nabla P_{w}, \label{eq:single_phase_eq3}
\end{align}
where
we assume $q$ is a function independent of time such that we can
compute the velocity $u_w$ in the initial time step and use the velocity
to construct the basis function for saturation $S_w$.

For solving single phase flow problem, we can separate the computation process into two parts. Since Equations (\ref{eq:single_phase_eq2}) and  (\ref{eq:single_phase_eq3}) are independent of saturation and the given source term is independent of time, we know that the velocity is independent of time. Therefore, we can decouple system of equation into two parts which are the Equations (\ref{eq:single_phase_eq2})-(\ref{eq:single_phase_eq3}) and the Equation (\ref{eq:single_phase_eq1}). The first part of the computation is computing the velocity and pressure by solving Equations (\ref{eq:single_phase_eq2}) and  (\ref{eq:single_phase_eq3}) by Mixed Generalized Multiscale Finite Element Methods(MGMsFEM) \cite{MixedGMsFEM}.
The second part of the computation is computing the numerical solution for the saturation by solving Equation (\ref{eq:single_phase_eq1}). Given a numerical velocity $u_{ms}$, Equation (\ref{eq:single_phase_eq1}) is a standard linear transport equation. We can use the method discussed in Section \ref{sec:lin-tran} to construct the multiscale basis function and then compute the numerical solution for saturation. Next, we will discuss the computation of the saturation.

We will consider a multi-continuum version of the method. Assume that an approximation of the velocity $u_{ms}$ has been computed by the MGMsFEM.
We derive the upscaled system following the ideas in Section~\ref{sec:overview}. For each coarse element $K_i$ and a continuum $j$ within $K_i$,
we solve the following system in an oversampled region $K_i^+$:
\begin{equation}
\label{eq:lin-trans-2}
\nabla \cdot (u_{ms} \psi^j_i) + \mu = 0
\end{equation}
to obtain a basis function $\psi^j_i$. The above system is equipped with the constraints
\begin{equation}
\label{eq:lin-trans-1}
\int_{K_m^{(l)}} \psi^j_i = \delta_{mj}\delta_{il}.
\end{equation}
We remark that $\mu$ is a piecewise constant function, and plays the role of Lagrange multiplier.
We notice that the upscaled system has the form (\ref{eq:nlmc-linear}).

\begin{figure}[!ht]
\centering
\includegraphics[scale=0.4]{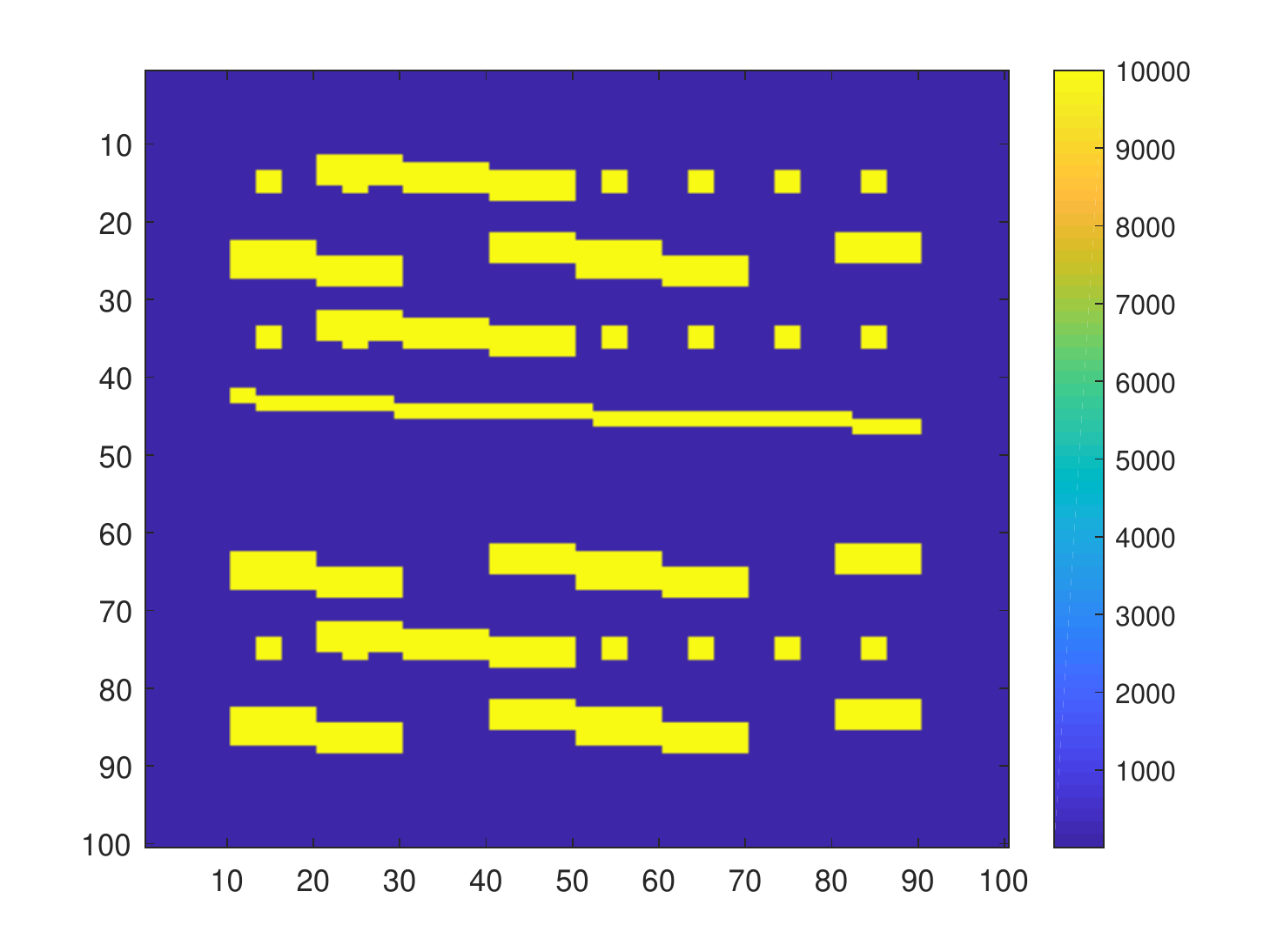}
\caption{permeability field $\kappa$}
\label{fig:permeability_case1}
\end{figure}
\begin{figure}[!ht]
\centering
\includegraphics[scale=0.35]{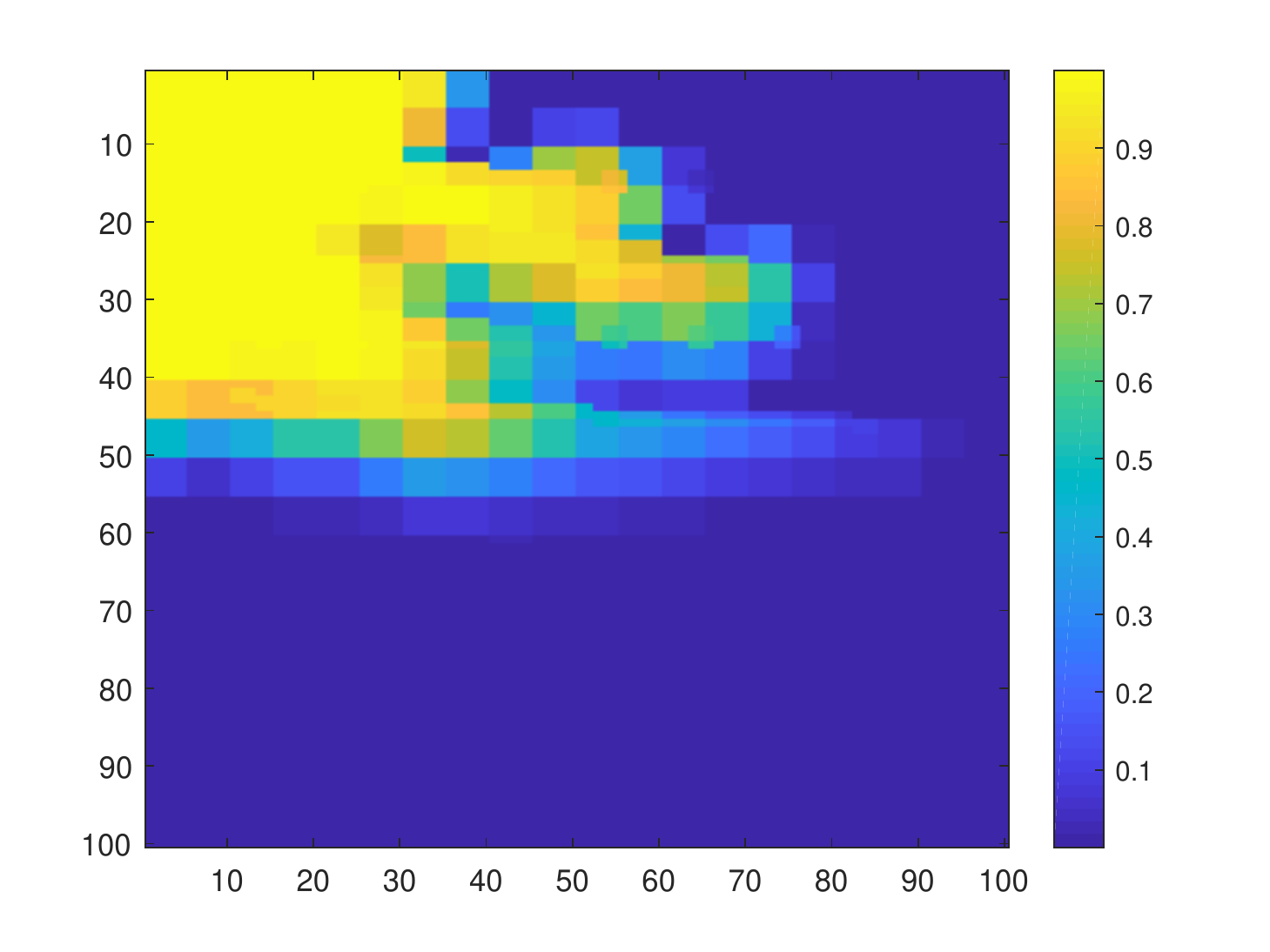} \includegraphics[scale=0.35]{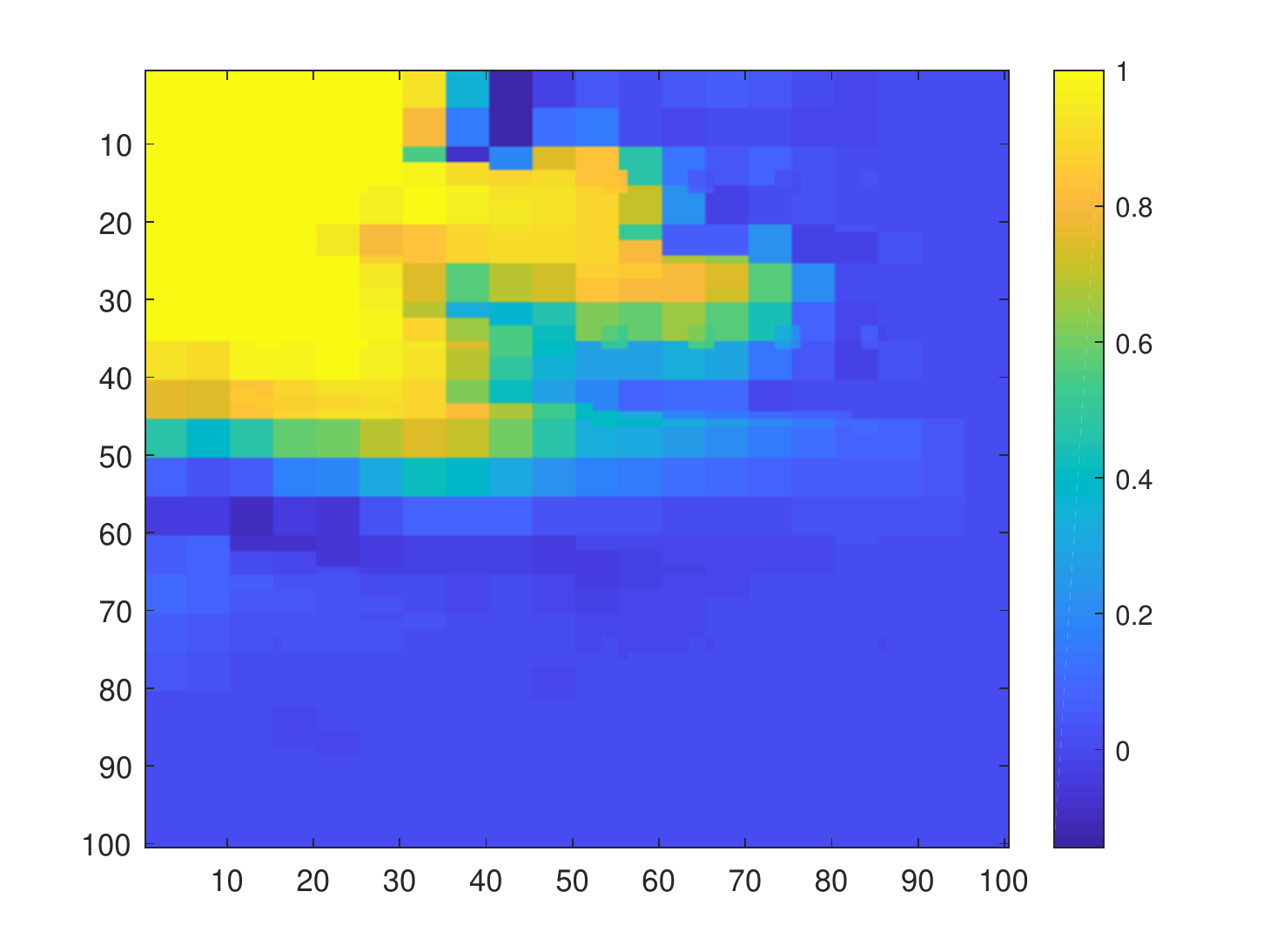}
\includegraphics[scale=0.35]{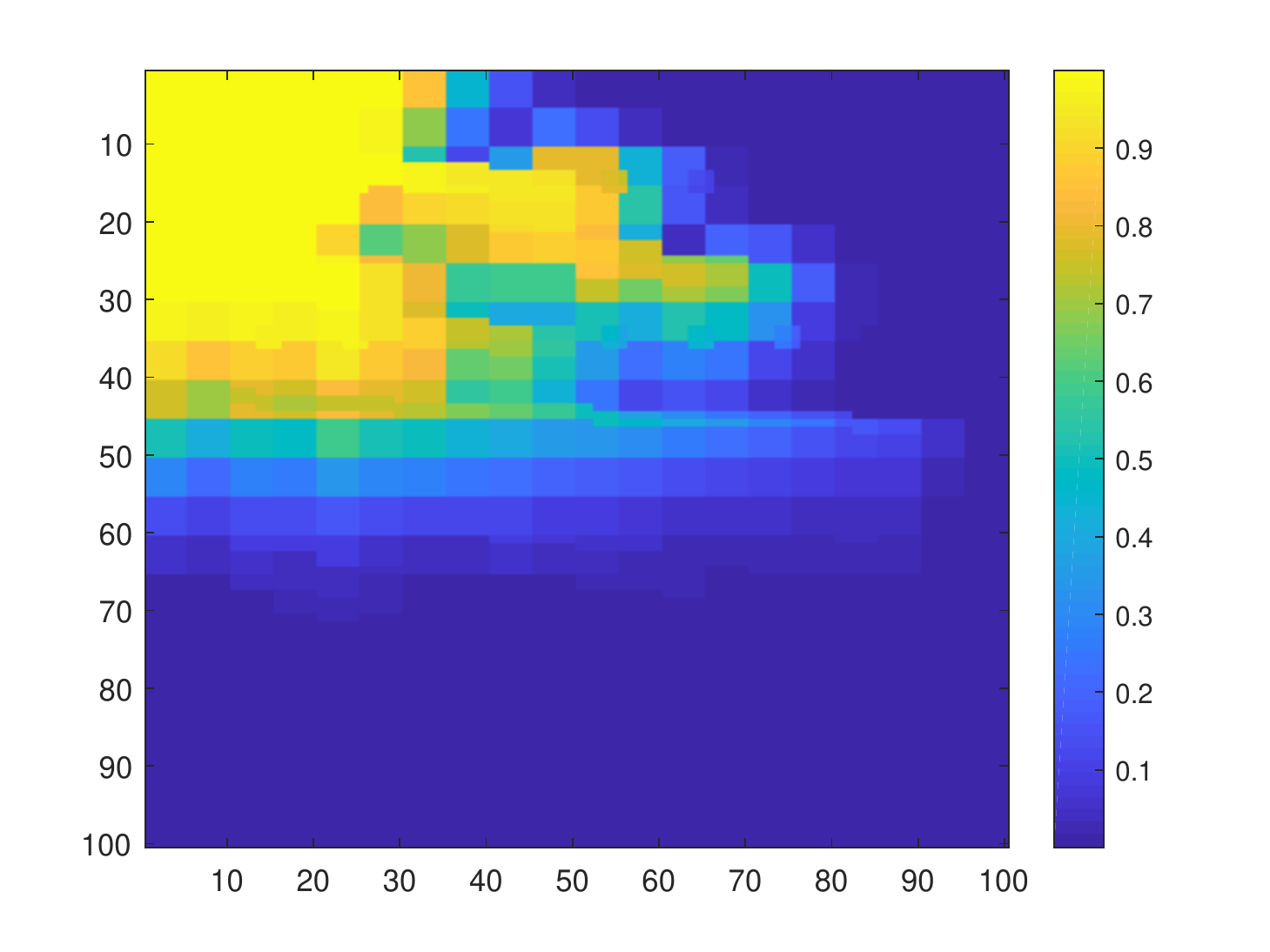}
\caption{Snapshots of the solution at $T=10$. Left: fine solution, Middle:
upscaled solution, Right: finite volume solution.}
\label{fig:ex13}
\end{figure}
\begin{figure}[!ht]
\centering
\includegraphics[scale=0.35]{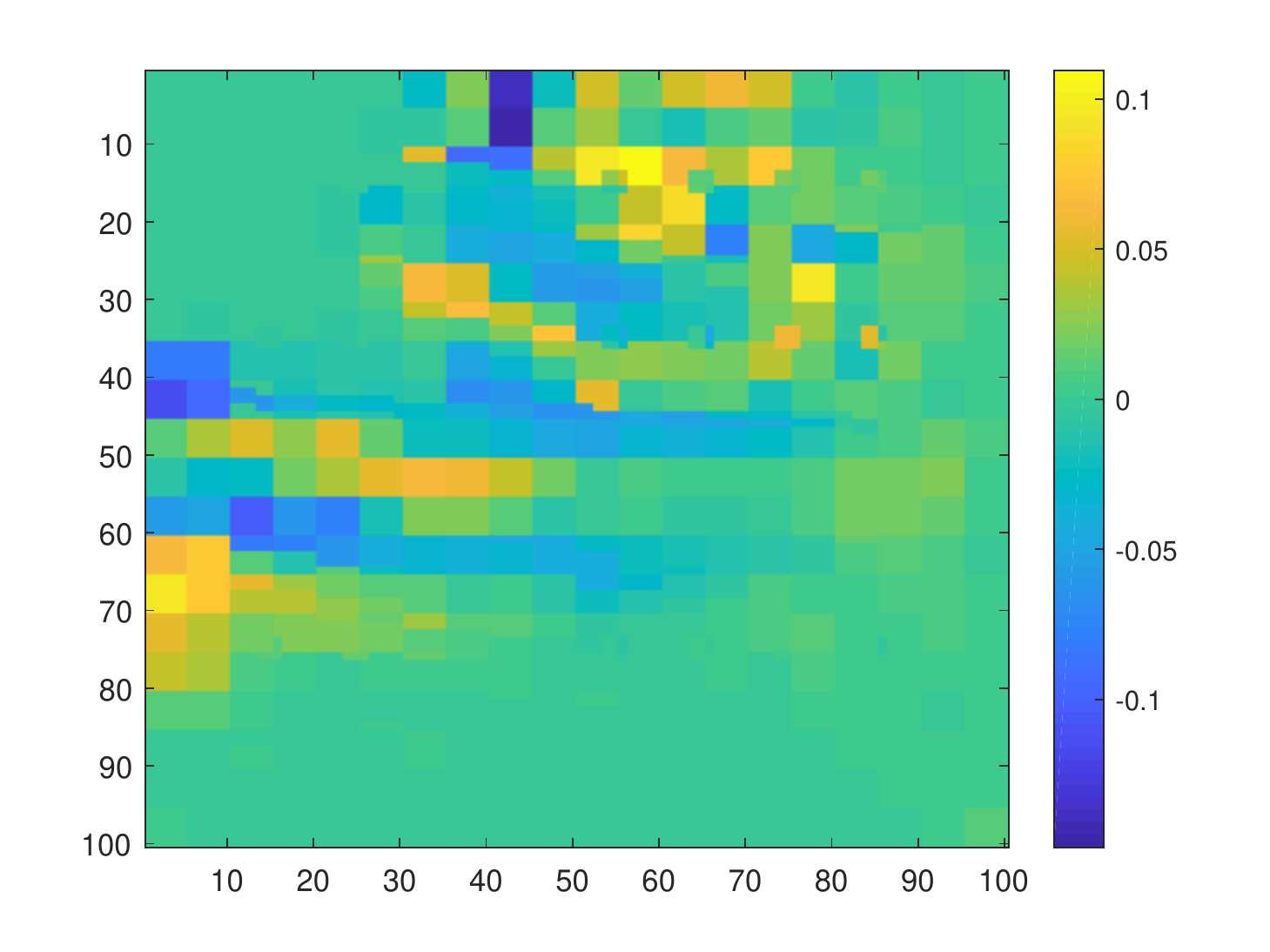} \includegraphics[scale=0.35]{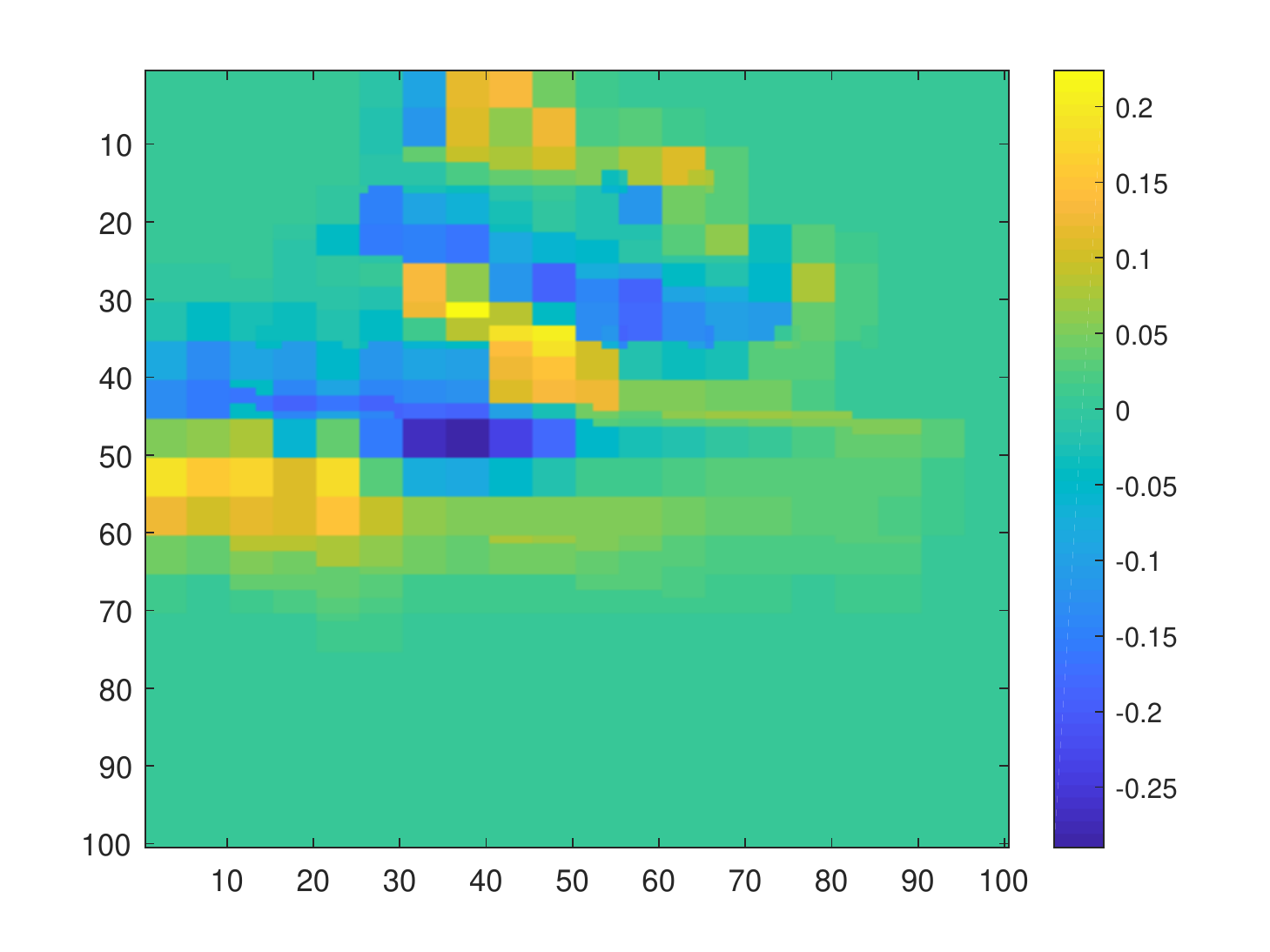}

\caption{Error comparison to fine solution at $T=10$. Left: upscaled solution,
Right: finite volume solution.}
\label{fig:ex14}
\end{figure}

Now, we present some numerical examples.
We take $\Omega=[0,1]^{2}$, and
the permeability field $\kappa$ is shown in Figure \ref{fig:permeability_case1}.
The coarse grid size $H=1/20$. We will use a dual continuum model
to solve this problem.
For each coarse element $K_j$,
we define two continua by $K_{j}^{(1)}=\{x\in K_{j}|\kappa(x)>10^{3}\}$
and $K_{j}^{(2)}=\{x\in K_{j}|\kappa(x) \leq 10^{3}\}$.
In Table \ref{tab:single_error_case1}, we present a error comparison at times $T=5$ and $T=10$
between our method and the finite volume method for the saturation
equation (\ref{eq:single_phase_eq3}) with piecewise constant test functions in each continuum.
From the table, we observe the performance of our scheme for various choices of oversampling layers.
For all these cases, we observe that our scheme performs better than the usual (dual continuum) finite volume scheme
whose error is $13.10\%$ and $14.04\%$ at the observation times $T=5$ and $T=10$ respectively.
We notice that the relative error in the velocity $u_{ms}$ is $5.51\%$.
In Figure \ref{fig:ex13} - \ref{fig:ex14}, we present the comparisons
between the the averaged fine-grid solution and the solution obtained
our proposed method. From these figures, we observe that our proposed
method can capture the solution features accurately.

\begin{table}[!ht]
\centering
\begin{tabular}{|c|c|c|}
\hline
\#layer \textbackslash{} time & $T=5$ & $T=10$\tabularnewline
\hline
\hline
4 & 12.83\% & 14.57\%\tabularnewline
\hline
6 & 7.33\% & 6.93\%\tabularnewline
\hline
8 & 5.89\% & 5.51\%\tabularnewline
\hline
$\infty$ & 5.69\% & 5.24\%\tabularnewline
\hline
\end{tabular}
\caption{$L^2$ relative errors for our upscaling method. The relative error for dual continuum finite volume scheme on the same grid is $13.10\%$ and $14.04\%$
for $T=5$ and $T=10$ respectively. In this case, the relative error in velocity $u_{ms}$ is $5.51\%$.}
\label{tab:single_error_case1}
\end{table}

%
%

\subsection{Two phase flow}

In this section, we apply our upscaling method to two phase flow problem which is
described as follows:
\begin{align}
\partial_{t}S_{w}+\nabla\cdot(u_{w}) & =q_{w}, \label{eq:two_phase_eq1}\\
-\nabla\cdot(\kappa\lambda_{t}(S_{w})\nabla P_{w}) & =q, \label{eq:two_phase_eq2}\\
u_{w} & =-\kappa\lambda_{w}(S_{w})\nabla P_{w}, \label{eq:two_phase_eq3}\\
\lambda_{t} & =\lambda_{w}(S_{w})+\lambda_{o}(S_{w}).
\end{align}
In the process of basis construction, we will use a similar approach as solving single phase flow problem. We notice that the exact velocity depends on the saturation and hence depends on time. To reduce the computational burden, we use the initial total mobility $\lambda_t(S_{w}(0))$ as the reference mobility to construct the multiscale basis functions. That is, we will construct the multiscale basis functions for approximating the exact velocity field based on the following equations,
\begin{equation}
\label{eq:two_phase_eq5}
-\nabla\cdot(u_t) =q,
\end{equation}
\begin{equation}
\label{eq:two_phase_eq6}
u_{t} =-\kappa\lambda_t(S_w(0,\cdot))\nabla P.
\end{equation}
Next, using the total velocity $u_t$, we will construct the basis functions for approximating the exact saturation based on the following transport equation,
 \begin{align}
\partial_{t}S_{w}+\nabla\cdot(\cfrac{u_{t}}{\lambda_t(S_w(0,\cdot))}S_w) & =q_{w}, \label{eq:two_phase_eq7}.
\end{align}
Hence, the construction for the basis functions is separated into two parts. The first part is computing an approximating velocity $u^{(0)}_{t,ms}$ and pressure $P_{ms}$ by solving Equations (\ref{eq:two_phase_eq5}) and  (\ref{eq:two_phase_eq6}) by Mixed Generalized Multiscale Finite Element Methods (MGMsFEM). The multiscale basis functions constructed during the processing are used to span the finite element space $V_{ms}$, which is used to approximate the exact velocity.
The second part is constructing the saturation basis functions by the method discussed in Section \ref{sec:concept} with a fixed velocity $u^{(0)}_{t,ms}$.
This part is also similar to the single phase flow.

To compute the coarse grid solution, we will use IMplicit
Pressure Explicit Saturation (IMPES) scheme. Given the saturation $S^{(n)}_{ms}$ and the total velocity $u^{(n)}_{t,ms}$ at the previous time step, we will compute the total velocity at the next time step $u^{(n+1)}_{t,ms}$ by solving
 \begin{align*}
\int_{\Omega} \nabla\cdot u^{(n+1)}_{t,ms}p & =\int_{\Omega} qp,\quad \;\forall p\in Q_H\\
\int_{\Omega} \kappa^{-1}\lambda^{-1}_t(S^{(n)}_w(,\cdot)) u^{(n+1)}_{t,ms}\cdot v & = \int_{\Omega} P^{(n+1)}_{ms},\nabla\cdot v,\quad\;\forall v\in V_{ms}
\end{align*}
where $Q_H$ is the space for pressure.
Then, we will use the velocity $u^{(n+1)}_{t,ms}$ to compute the saturation at the next time step
as before (similar to (\ref{eq:lin-trans-2}) and (\ref{eq:lin-trans-1})).



\begin{figure}[!ht]
\centering
\includegraphics[scale=0.4]{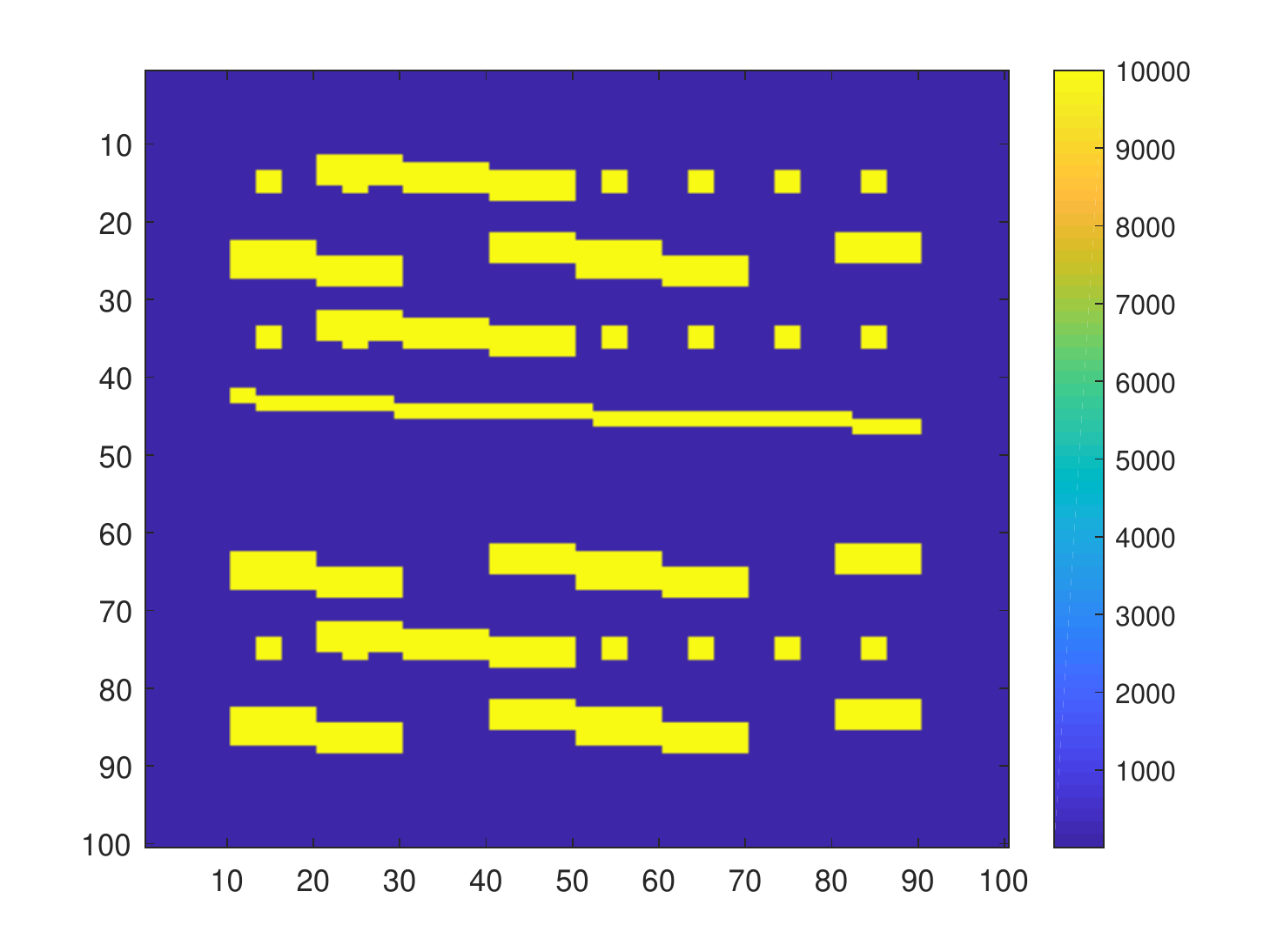}
\caption{permeability field $\kappa$}
\label{fig:permeability_case1-1}
\end{figure}

We next present some numerical examples.
We take $\Omega=[0,1]^{2}$.
The permeability field $\kappa$ is shown in Figure \ref{fig:permeability_case1-1}, and we use
\[
\lambda_{w}=\Big(\cfrac{S_{w}-S_{wc}}{(1-S_{wc}-S_{or})}\Big)^{2}\mu_{w}^{-1}, \quad \lambda_{o}=\Big(\cfrac{1-S_{w}-S_{or}}{1-S_{wc}-S_{or}}\Big)^{2}\mu_{o}^{-1}
\]
with $\mu_{w}=1,$ $\mu_{o}=1$, $S_{wc}=0.2$ and $S_{or}=0.2$.
The coarse grid size $H=1/20$. We are using the dual continuum model
to solve this problem. For each coarse element $K_j$, we define two continua by $K_{j}^{(1)}=\{x\in K_{j}|\kappa(x)>10^{3}\}$
and $K_{j}^{(2)}=\{x\in K_{j}|\kappa(x) \leq 10^{3}\}$.
In Table \ref{tab:single_error_case1-1}, we present a error comparison
between our method and the finite volume method for the saturation
equation with piecewise constant test functions in each continuum.
For this experiment, we observed that the upscaling method works in general better than the finite volume method as before.

\begin{table}[!ht]
\centering
\begin{tabular}{|c|c|c|}
\hline
\#layer \textbackslash{} time & $T=10$ & $T=20$\tabularnewline
\hline
\hline
4 & 15.82\% & 14.94\%\tabularnewline
\hline
6 & 7.03\% & 6.63\%\tabularnewline
\hline
8 & 4.11\% & 4.33\%\tabularnewline
\hline
$\infty$ & 3.52\% & 4.11\%\tabularnewline
\hline
\end{tabular}
\caption{$L^2$ relative errors for our upscaling method. The relative error for dual continuum finite volume scheme on the same grid is $10.01\%$ and $14.13\%$
for $T=10$ and $T=20$ respectively.}
\label{tab:single_error_case1-1}
\end{table}

\begin{figure}[!ht]
\centering

\includegraphics[scale=0.4]{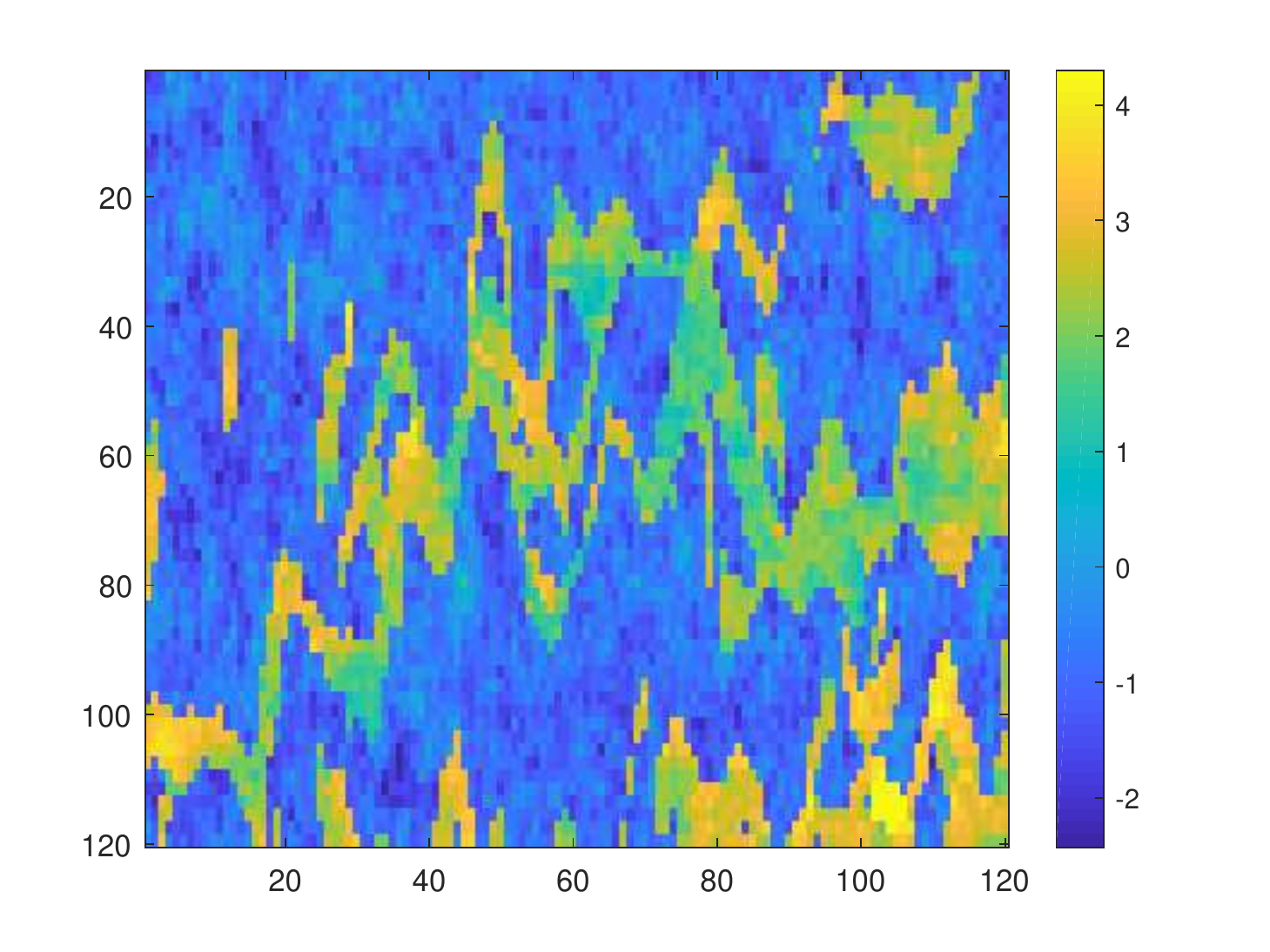}

\caption{A SPE permeability field $\log_{10}(\kappa)$ (in log scale).}
\label{fig:permeability_case2-1}
\end{figure}

\begin{figure}[!ht]
\centering
\includegraphics[scale=0.35]{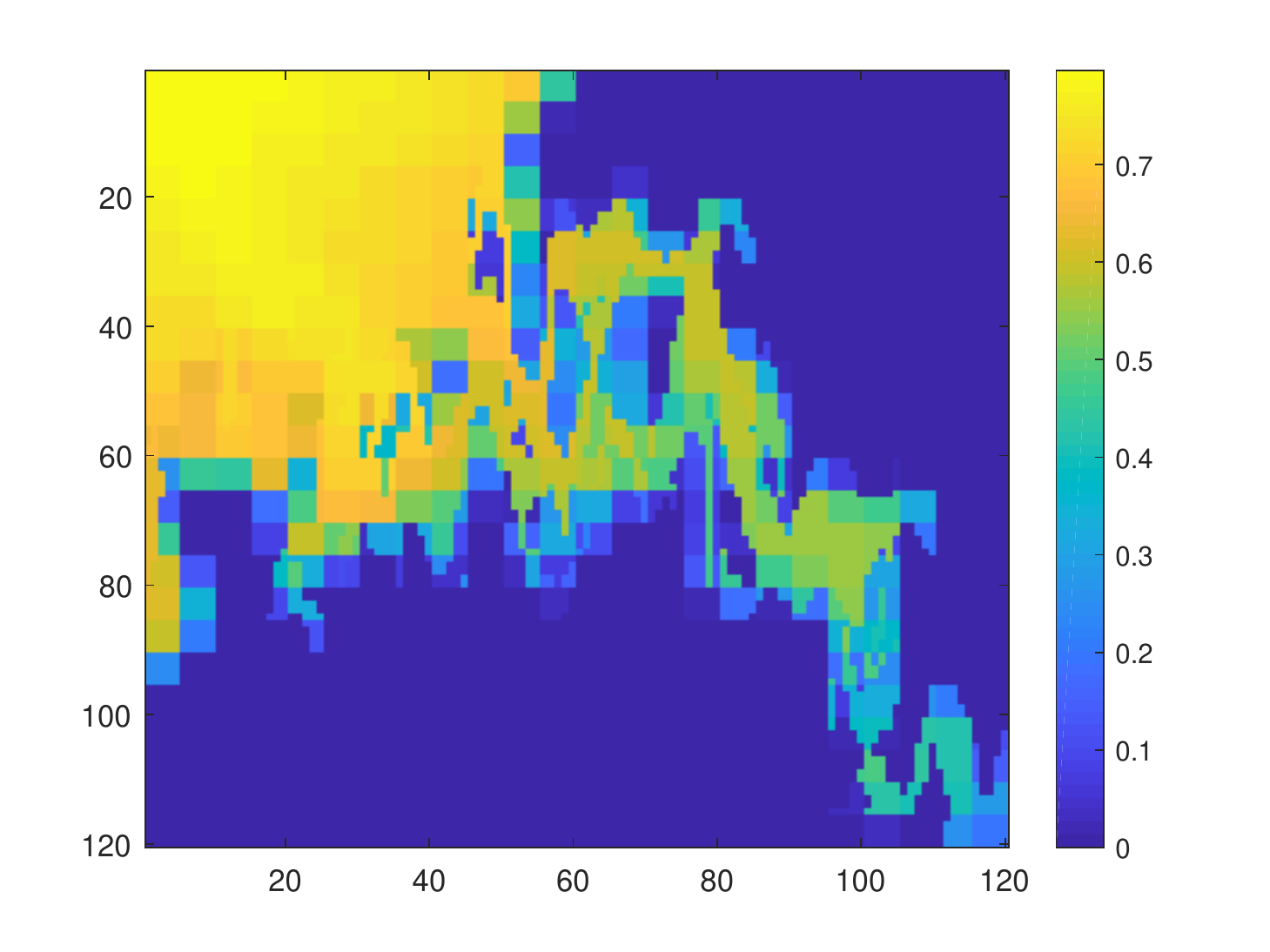}
\includegraphics[scale=0.35]{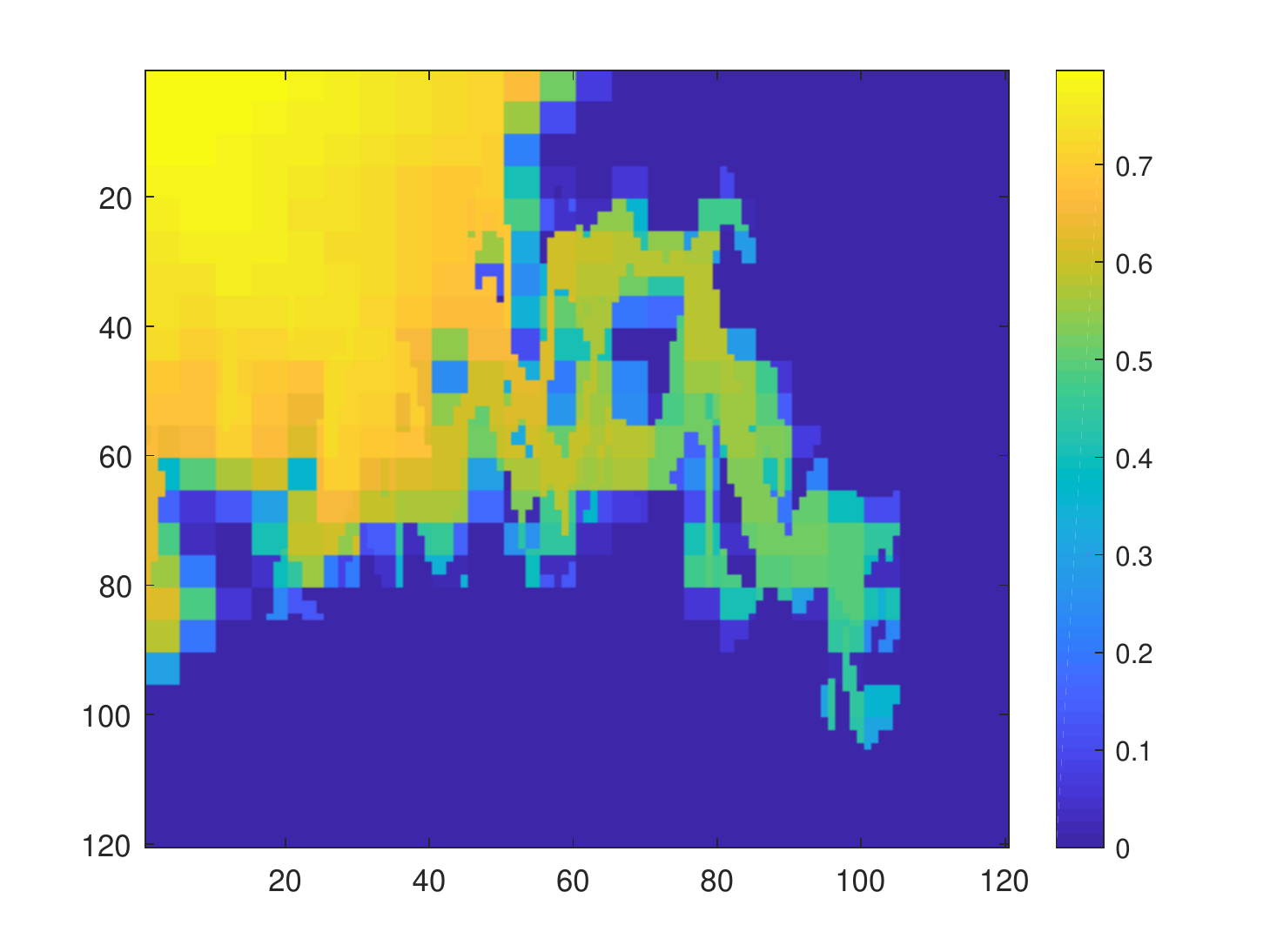}
\includegraphics[scale=0.35]{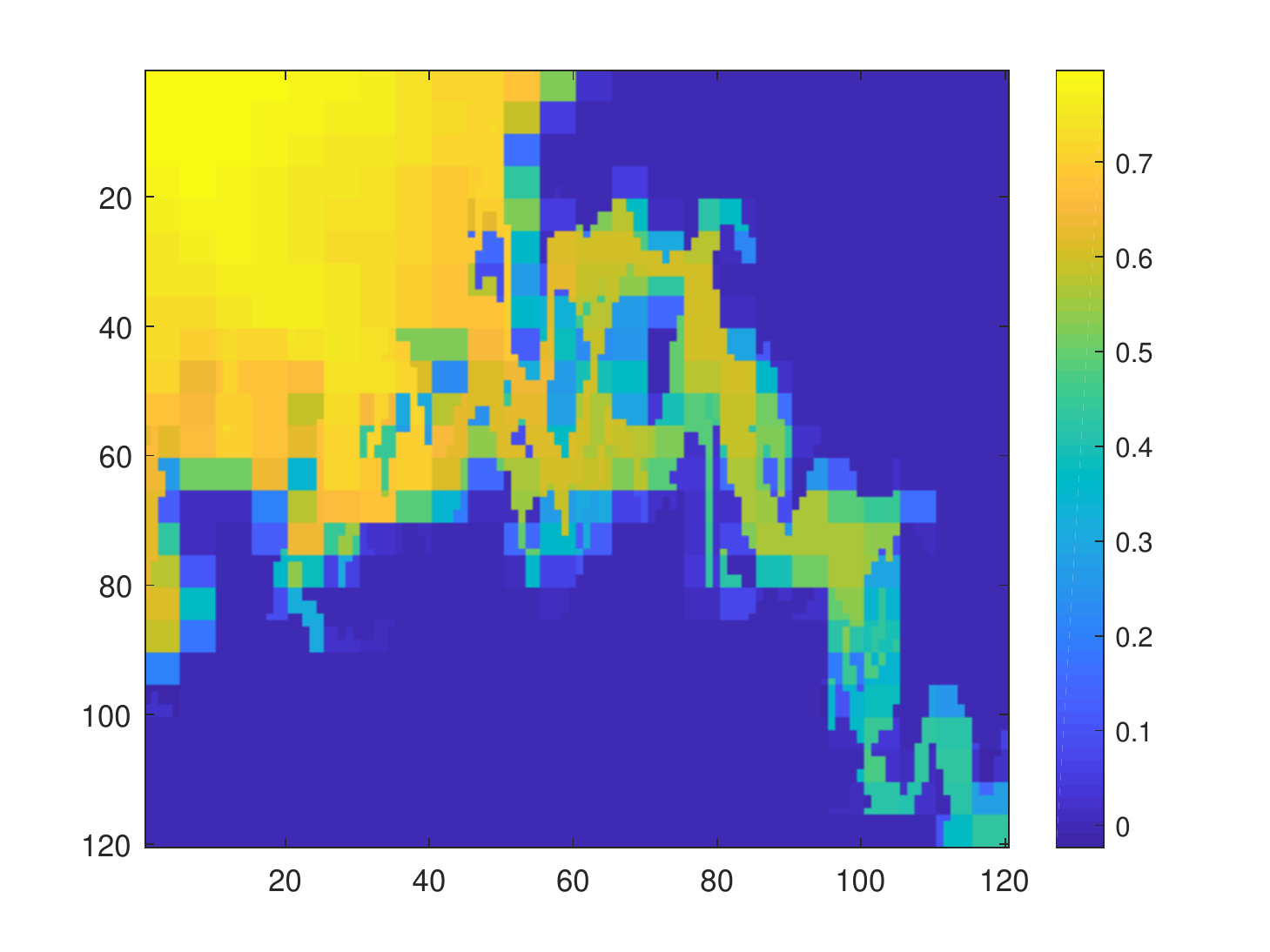}
\caption{Snapshots of the solution at $T=10$. Left: averaged fine solution. Middle: finite volume.
Right: our upscaled solution.}
\label{fig:2p-snap}
\end{figure}

Next we consider another test case with the medium $\kappa$ shown in Figure~\ref{fig:permeability_case2-1},
which is a SPE benchmark test case.
The coarse grid size is chosen as $H=1/24$. We are again using the dual continuum model
to solve this problem. For each coarse element $K_j$, we define
$K_{j}^{(1)}=\{x\in K_{j}|\log_{10}(\kappa(x))>0.8\}$ and $K_{j}^{(2)}=\{x\in K_{j}|\log_{10}(\kappa(x)) \leq 0.8\}$.
In Figure~\ref{fig:2p-snap}, we present the snapshots of the averaged fine solution, finite volume solution
and our upscaled solution at the observation time $T=10$.
We observe that our method is able to capture the dynamics of the solution, while the finite volume
method does not produce a good solution.
In Table~\ref{tab:single_error_case2-1}, we present the errors of our approximation and observe again good accuracy.

\begin{table}[!ht]
\centering

\begin{tabular}{|c|c|c|}
\hline
\#layer \textbackslash{} time & $T=5$ & $T=10$\tabularnewline
\hline
6 & 7.19\% & 11.12\%\tabularnewline
\hline
$\infty$ & 6.70\% & 10.16\%\tabularnewline
\hline
\end{tabular}

\caption{$L^2$ relative errors for our upscaling method. The relative error for dual continuum finite volume scheme on the same grid is $14.75\%$ and $19.74\%$
for $T=5$ and $T=10$ respectively.}
\label{tab:single_error_case2-1}
\end{table}

\section{Nonlinear Approach}
 \label{sec:nonlinear}

In this section, we will present a numerical test to show the performance of the nonlinear approach (c.f. Section \ref{sec:concept}).
To to do, we consider the following
single phase flow problem
\begin{align*}
\partial_{t}S_{w}+\nabla\cdot(u_{w}\lambda (S_{w})) & =q_{w},\\
-\nabla\cdot(\kappa\nabla P_{w}) & =q,\\
u_{w} & =-\kappa\nabla P_{w},
\end{align*}
where
\[
\lambda (S) = S^{\beta}.
\]
We assume $q$ is a function independent of time such that we can
compute the velocity in the initial time step and use the velocity
to construct the basis function for saturation.

The derivation of the upscaled system follows the ideas in Section~\ref{sec:non-nlmc}.
Assume that an approximation of the velocity $u_{ms}$ has been computed.
Let $\{ S^{n,j}_i\}$ be a set of upscaled values for the continuum $j$ in the coarse region $K_i$ at the time $t_n$.
For each coarse element $K_i$,
we solve the following system in an oversampled region $K_i^+$:
\begin{equation}
\label{eq:non-trans}
\nabla \cdot (u_{ms}  \lambda(N_i) ) + \mu = 0
\end{equation}
to obtain a local downscale function $N_i$. The above system is equipped with the constraints
\begin{equation}
\frac{1}{| K_i^{(j)}|} \int_{K_i^{(j)}} N_i = S^{n,j}_i.
\end{equation}
We remark that $\mu$ is a piecewise constant function, and plays the role of Lagrange multiplier.
Then we define a global downscale field by
$$
S^n_h = \sum_i N_i \chi^{K_i^+},
$$
where $\{ \chi^{K_i^+}\}$ is a set of partition of unity functions corresponding to the sets $K_i^+$.  Then we apply a finite volume scheme
to the saturation equation as follows
\begin{equation}
\label{eq:nonlinear-1p}
S^{n+1,j}_i = S^{n,j}_i - \Delta t \int_{\partial K_i^{(j)}} \lambda(S^n_h) \, u_{ms}\cdot n + \Delta t \int_{K_i^{(j)}} q_w.
\end{equation}
We notice that the upscaled system (\ref{eq:nonlinear-1p}) has the form (\ref{eq:nlnlmc}).

\begin{figure}[!ht]
\centering
\includegraphics[scale=0.4]{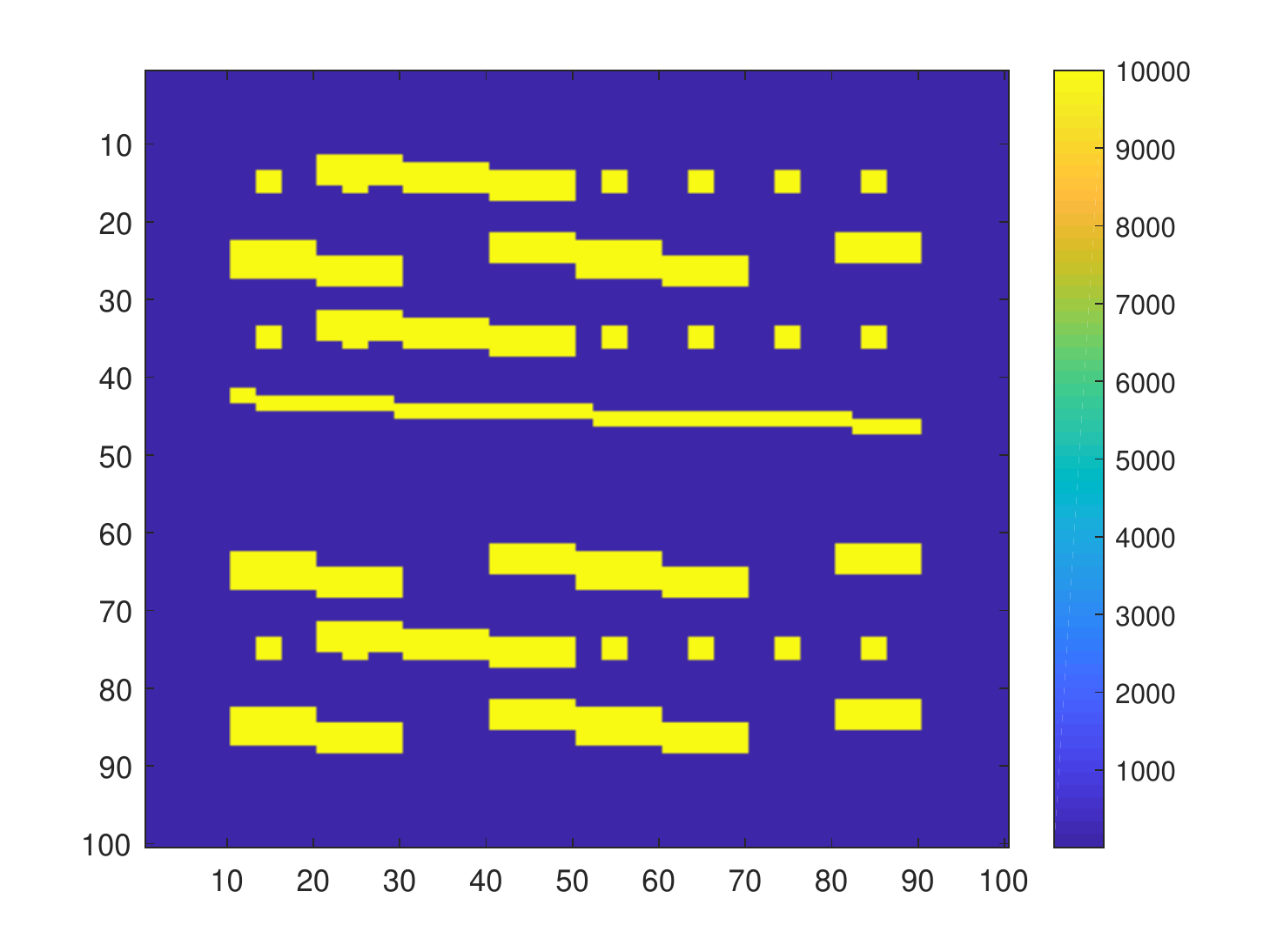}
\caption{permeability field $\kappa$}
\label{fig:permeability_case1-2}
\end{figure}

We next present some numerical results.
We take $\Omega=[0,1]^{2}$, and
the permeability field $\kappa$ is shown in Figure \ref{fig:permeability_case1-2}.
The coarse grid size $H=1/10$. We are using a dual continuum model
to solve this problem. For each $K_j$, we define $K_{j}^{(1)}=\{x\in K_{j}|\kappa(x)>10^{3}\}$
and $K_{j}^{(2)}=\{x\in K_{j}|\kappa(x) \leq 10^{3}\}$.
In Table \ref{tab:single_error_case1-2}, we present a error comparison
between our method and the finite volume method for the saturation
equation with piecewise constant test functions in each continua.
We consider three choices of $\beta$, and present the relative errors at times $T=5$ and $T=10$.
From the results, we observe that our method is able to compute more accurate numerical solutions.

\begin{table}[!ht]
\centering
\begin{tabular}{|c|c|c|}
\hline
$\beta$ \textbackslash{} time & $T=5$ & $T=10$\tabularnewline
\hline
\hline
2 & 9.91\% & 11.85\%\tabularnewline
\hline
3 & 7.54\% & 12.60\%\tabularnewline
\hline
5 & 10.07\% & 14.89\%\tabularnewline
\hline
\end{tabular}
\begin{tabular}{|c|c|c|}
\hline
$\beta$ \textbackslash{} time & $T=5$ & $T=10$\tabularnewline
\hline
\hline
2 & 28.70\% & 23.61\%\tabularnewline
\hline
3 & 31.32\% & 26.27\%\tabularnewline
\hline
5 & 32.98\% & 29.19\%\tabularnewline
\hline
\end{tabular}
\caption{$L^2$ norm relative errors. Left: Nonlinear upscaling. Right: Finite volume method.}
\label{tab:single_error_case1-2}
\end{table}

We remark that our proposed upscaling method can be applied to other nonlinear systems,
which will be investigated in more detail in a forthcoming paper.

\section{Conclusions}

In this paper, we study and develop
a general framework for coarse-grid modeling
of nonlinear PDEs. The fine-grid problem is described by nonlinear
PDEs including two-phase flow and transport model.
The main concept of our upscaled model is three fold,
(1) We determine macroscopic quantities in each coarse-grid block,
which can vary among the coarse-grid blocks.
(2) We define appropriate local problems with constraints
 in oversampled regions
which are solvable. The constraints are given by macroscopic quantities.
These local problems allow downscaling from macroscopic quantities to
each coarse-grid block.
(3)  We formulate the coarse-grid problem such that the downscaled
solution over the entire domain solves the global problem in a weak
sense.

The formulation of local problems allow easily defining macroscopic
fluxes. The local problems are defined via forcing terms. We
also discuss the use of multiscale basis function in a linear fashion
for solving nonlinear PDEs.
Our main example includes two-phase flow and transport and its simplifications.
The model problem is described by a coupled PDEs. We consider various permeability fields and show that one can achieve an accurate approximation of the solution.

\section*{Acknowledgements}

EC's work is partially supported by Hong Kong RGC General Research Fund (Project 14304217)
and CUHK Direct Grant for Research 2017-18. EC would like to thank the support of the ESI
for attending the program Numerical Analysis of Complex PDE Models in the Sciences.

\bibliographystyle{plain}
\bibliography{references,references1,references2}

\begin{thebibliography}{10}

\bibitem{abdul_yun}
Assyr Abdulle and Yun Bai.
\newblock Adaptive reduced basis finite element heterogeneous multiscale
  method.
\newblock {\em Comput. Methods Appl. Mech. Engrg.}, 257:203--220, 2013.

\bibitem{ab05}
G.~Allaire and R.~Brizzi.
\newblock A multiscale finite element method for numerical homogenization.
\newblock {\em SIAM J. Multiscale Modeling and Simulation}, 4(3):790--812,
  2005.

\bibitem{arbogast02}
T.~Arbogast.
\newblock Implementation of a locally conservative numerical subgrid upscaling
  scheme for two-phase {D}arcy flow.
\newblock {\em Comput. Geosci}, 6:453--481, 2002.

\bibitem{apwy07}
T.~Arbogast, G.~Pencheva, M.F. Wheeler, and I.~Yotov.
\newblock A multiscale mortar mixed finite element method.
\newblock {\em SIAM J. Multiscale Modeling and Simulation}, 6(1):319--346,
  2007.

\bibitem{Arbogast_PWY_07}
T.~Arbogast, G.~Pencheva, M.F. Wheeler, and I.~Yotov.
\newblock A multiscale mortar mixed finite element method.
\newblock {\em Multiscale Model. Simul.}, 6(1):319--346, 2007.

\bibitem{barenblatt1960basic}
GI~Barenblatt, Iu~P Zheltov, and IN~Kochina.
\newblock Basic concepts in the theory of seepage of homogeneous liquids in
  fissured rocks [strata].
\newblock {\em Journal of applied mathematics and mechanics}, 24(5):1286--1303,
  1960.

\bibitem{BT}
J.W. Barker and S.~Thibeau.
\newblock A critical review of the use of pseudorelative permeabilities for
  upscaling.
\newblock {\em SPE Reservoir Eng.}, 12:138--143, 1997.

\bibitem{brown2014multiscale}
Donald~L Brown and Daniel Peterseim.
\newblock A multiscale method for porous microstructures.
\newblock {\em arXiv preprint arXiv:1411.1944}, 2014.

\bibitem{cdgw03}
Y.~Chen, L.~Durlofsky, M.~Gerritsen, and X.~Wen.
\newblock A coupled local-global upscaling approach for simulating flow in
  highly heterogeneous formations.
\newblock {\em Advances in Water Resources}, 26:1041--1060, 2003.

\bibitem{ElasticGMsFEM}
E.~Chung, Y.~Efendiev, and S.~Fu.
\newblock Generalized multiscale finite element method for elasticity
  equations.
\newblock {\em International Journal on Geomathematics}, 5(2):225--254, 2014.

\bibitem{MixedGMsFEM}
E.~Chung, Y.~Efendiev, and C.~Lee.
\newblock Mixed generalized multiscale finite element methods and applications.
\newblock {\em SIAM Multicale Model. Simul.}, 13:338--366, 2014.

\bibitem{WaveGMsFEM}
E.~Chung, Y.~Efendiev, and W.~T. Leung.
\newblock Generalized multiscale finite element method for wave propagation in
  heterogeneous media.
\newblock {\em SIAM Multicale Model. Simul.}, 12:1691--1721, 2014.

\bibitem{MsDG}
E.~Chung and W.~T. Leung.
\newblock A sub-grid structure enhanced discontinuous galerkin method for
  multiscale diffusion and convection-diffusion problems.
\newblock {\em Communications in Computational Physics}, 14:370--392, 2013.

\bibitem{OnlineStokes}
E.~T. Chung, Y.~Efendiev, W.T. Leung, M.~Vasilyeva, and Y.~Wang.
\newblock Online adaptive local multiscale model reduction for heterogeneous
  problems in perforated domains.
\newblock {\em Applicable Analysis}, 96(12):2002--2031, 2017.

\bibitem{AdaptiveGMsFEM}
E.~T. Chung, Y.~Efendiev, and G.~Li.
\newblock An adaptive {GM}s{FEM} for high contrast flow problems.
\newblock {\em J. Comput. Phys.}, 273:54--76, 2014.

\bibitem{chung2018constraintmixed}
Eric Chung, Yalchin Efendiev, and Wing~Tat Leung.
\newblock Constraint energy minimizing generalized multiscale finite element
  method in the mixed formulation.
\newblock {\em Computational Geosciences}, 22(3):677--693, 2018.

\bibitem{chung2017DGstokes}
Eric Chung, Maria Vasilyeva, and Yating Wang.
\newblock A conservative local multiscale model reduction technique for stokes
  flows in heterogeneous perforated domains.
\newblock {\em Journal of Computational and Applied Mathematics}, 321:389--405,
  2017.

\bibitem{NLMC}
Eric~T Chung, Efendiev, Wing~Tat Leung, Maria Vasilyeva, and Yating Wang.
\newblock Non-local multi-continua upscaling for flows in heterogeneous
  fractured media.
\newblock {\em arXiv preprint arXiv:1708.08379}, 2018.

\bibitem{chung2018constraint}
Eric~T Chung, Yalchin Efendiev, and Wing~Tat Leung.
\newblock Constraint energy minimizing generalized multiscale finite element
  method.
\newblock {\em Computer Methods in Applied Mechanics and Engineering},
  339:298--319, 2018.

\bibitem{chung2018fast}
Eric~T Chung, Yalchin Efendiev, and Wing~Tat Leung.
\newblock Fast online generalized multiscale finite element method using
  constraint energy minimization.
\newblock {\em Journal of Computational Physics}, 355:450--463, 2018.

\bibitem{ohl12}
Martin Drohmann, Bernard Haasdonk, and Mario Ohlberger.
\newblock Reduced basis approximation for nonlinear parametrized evolution
  equations based on empirical operator interpolation.
\newblock {\em SIAM J. Sci. Comput.}, 34(2):A937--A969, 2012.

\bibitem{dur91}
L.J. Durlofsky.
\newblock Numerical calculation of equivalent grid block permeability tensors
  for heterogeneous porous media.
\newblock {\em Water Resour. Res.}, 27:699--708, 1991.

\bibitem{ee03}
W.~E and B.~Engquist.
\newblock Heterogeneous multiscale methods.
\newblock {\em Comm. Math. Sci.}, 1(1):87--132, 2003.

\bibitem{ed01}
Y.~Efendiev and L.J. Durlofsky.
\newblock Numerical modeling of subgrid heterogeneity in two phase flow
  simulations.
\newblock {\em Water Resour. Res.}, 38(8):1128, 2002.

\bibitem{ed03}
Y.~Efendiev and L.J. Durlofsky.
\newblock A generalized convection-diffusion model for subgrid transport in
  porous media.
\newblock {\em SIAM J. Multiscale Modeling and Simulation}, 1(3):504--526,
  2003.

\bibitem{GMsFEM13}
Y.~Efendiev, J.~Galvis, and {T. Y.} Hou.
\newblock Generalized multiscale finite element methods (gmsfem).
\newblock {\em Journal of Computational Physics}, 251:116--135, 2013.

\bibitem{egw10}
Y.~Efendiev, J.~Galvis, and X.H. Wu.
\newblock Multiscale finite element methods for high-contrast problems using
  local spectral basis functions.
\newblock {\em Journal of Computational Physics}, 230:937--955, 2011.

\bibitem{eh09}
Y.~Efendiev and T.~Hou.
\newblock {\em {Multiscale Finite Element Methods: Theory and Applications}},
  volume~4 of {\em Surveys and Tutorials in the Applied Mathematical Sciences}.
\newblock Springer, New York, 2009.

\bibitem{ep03a}
Y.~Efendiev and A.~Pankov.
\newblock Numerical homogenization of monotone elliptic operators.
\newblock {\em SIAM J. Multiscale Modeling and Simulation}, 2(1):62--79, 2003.

\bibitem{ep03d}
Y.~Efendiev and A.~Pankov.
\newblock Numerical homogenization of nonlinear random parabolic operators.
\newblock {\em SIAM J. Multiscale Modeling and Simulation}, 2(2):237--268,
  2004.

\bibitem{ep03c}
Y.~Efendiev and A.~Pankov.
\newblock Homogenization of nonlinear random parabolic operators.
\newblock {\em Advances in Differential Equations}, 10(11):1235--1260, 2005.

\bibitem{fafalis2012capability}
DA~Fafalis, SP~Filopoulos, and GJ~Tsamasphyros.
\newblock On the capability of generalized continuum theories to capture
  dispersion characteristics at the atomic scale.
\newblock {\em European Journal of Mechanics-A/Solids}, 36:25--37, 2012.

\bibitem{fafalis2018computational}
Dimitrios Fafalis and Jacob Fish.
\newblock Computational continua for linear elastic heterogeneous solids on
  unstructured finite element meshes.
\newblock {\em International Journal for Numerical Methods in Engineering},
  115(4):501--530, 2018.

\bibitem{fish2013practical}
Jacob Fish.
\newblock {\em Practical multiscaling}.
\newblock John Wiley \& Sons, 2013.

\bibitem{fish2004space}
Jacob Fish and Wen Chen.
\newblock Space--time multiscale model for wave propagation in heterogeneous
  media.
\newblock {\em Computer Methods in applied mechanics and engineering},
  193(45):4837--4856, 2004.

\bibitem{fish2008mathematical}
Jacob Fish and Rong Fan.
\newblock Mathematical homogenization of nonperiodic heterogeneous media
  subjected to large deformation transient loading.
\newblock {\em International Journal for numerical methods in engineering},
  76(7):1044--1064, 2008.

\bibitem{fish2015computational}
Jacob Fish, Vasilina Filonova, and Dimitrios Fafalis.
\newblock Computational continua revisited.
\newblock {\em International Journal for Numerical Methods in Engineering},
  102(3-4):332--378, 2015.

\bibitem{fish2012reduced}
Jacob Fish, Vasilina Filonova, and Zheng Yuan.
\newblock Reduced order computational continua.
\newblock {\em Computer Methods in Applied Mechanics and Engineering},
  221:104--116, 2012.

\bibitem{fish2010computational}
Jacob Fish and Sergey Kuznetsov.
\newblock Computational continua.
\newblock {\em International Journal for Numerical Methods in Engineering},
  84(7):774--802, 2010.

\bibitem{fish1997computational}
Jacob Fish, Kamlun Shek, Muralidharan Pandheeradi, and Mark~S Shephard.
\newblock Computational plasticity for composite structures based on
  mathematical homogenization: Theory and practice.
\newblock {\em Computer Methods in Applied Mechanics and Engineering},
  148(1-2):53--73, 1997.

\bibitem{fish2005multiscale}
Jacob Fish and Zheng Yuan.
\newblock Multiscale enrichment based on partition of unity.
\newblock {\em International Journal for Numerical Methods in Engineering},
  62(10):1341--1359, 2005.

\bibitem{fish2007multiscale}
Jacob Fish and Zheng Yuan.
\newblock Multiscale enrichment based on partition of unity for nonperiodic
  fields and nonlinear problems.
\newblock {\em Computational Mechanics}, 40(2):249--259, 2007.

\bibitem{henning2009heterogeneous}
Patrick Henning and Mario Ohlberger.
\newblock The heterogeneous multiscale finite element method for elliptic
  homogenization problems in perforated domains.
\newblock {\em Numerische Mathematik}, 113(4):601--629, 2009.

\bibitem{hn00}
L.~Holden and B.F. Nielsen.
\newblock Global upscaling of permeability in heterogeneous reservoirs: the
  {\uppercase{o}}utput {\uppercase{l}}east {\uppercase{s}}quares
  ({\uppercase{ols}} method.
\newblock {\em Transport in Porous Media}, 40:115--143, 2000.

\bibitem{KB}
J.R. Kyte and D.W. Berry.
\newblock New pseudofunctions to control numerical dispersion.
\newblock {\em Society of Petroleum Engineers Journal}, 15(4):269--276, 1975.

\bibitem{lee2001hierarchical}
Seong~H Lee, MF~Lough, and CL~Jensen.
\newblock Hierarchical modeling of flow in naturally fractured formations with
  multiple length scales.
\newblock {\em Water resources research}, 37(3):443--455, 2001.

\bibitem{matache2002two}
Ana-Maria Matache and Christoph Schwab.
\newblock Two-scale fem for homogenization problems.
\newblock {\em ESAIM: Mathematical Modelling and Numerical Analysis},
  36(04):537--572, 2002.

\bibitem{oskay2007eigendeformation}
Caglar Oskay and Jacob Fish.
\newblock Eigendeformation-based reduced order homogenization for failure
  analysis of heterogeneous materials.
\newblock {\em Computer Methods in Applied Mechanics and Engineering},
  196(7):1216--1243, 2007.

\bibitem{oz07}
H.~Owhadi and L.~Zhang.
\newblock Metric-based upscaling.
\newblock {\em Comm. Pure. Appl. Math.}, 60:675--723, 2007.

\bibitem{panasenko2018multicontinuum}
GP~Panasenko.
\newblock Multicontinuum wave propagation in a laminated beam with contrasting
  stiffness and density of layers.
\newblock {\em Journal of Mathematical Sciences}, pages 1--13, 2018.

\bibitem{pankov97}
A.~Pankov.
\newblock {\em ${G}$-convergence and homogenization of nonlinear partial
  differential operators}.
\newblock Kluwer Academic Publishers, Dordrecht, 1997.

\bibitem{pwy02}
M.~Peszy\'nska, M.~Wheeler, and I.~Yotov.
\newblock Mortar upscaling for multiphase flow in porous media.
\newblock {\em Comput. Geosci.}, 6(1):73--100, 2002.

\bibitem{warren1963behavior}
JE~Warren, P~Jj Root, et~al.
\newblock The behavior of naturally fractured reservoirs.
\newblock 1963.

\bibitem{weh02}
X.H. Wu, Y.~Efendiev, and T.Y. Hou.
\newblock Analysis of upscaling absolute permeability.
\newblock {\em Discrete and Continuous Dynamical Systems, Series B.},
  2:158--204, 2002.

\bibitem{yuan2009multiple}
Zheng Yuan and Jacob Fish.
\newblock Multiple scale eigendeformation-based reduced order homogenization.
\newblock {\em Computer Methods in Applied Mechanics and Engineering},
  198(21-26):2016--2038, 2009.

\end{thebibliography}

\end{document}